\documentclass[11pt]{article}
\usepackage[centering,includeheadfoot,margin=2 cm]{geometry}
\usepackage{graphics,enumitem,epsfig}
\usepackage{amsfonts,amsmath,amssymb,euscript,color,stackrel}
\usepackage{hyperref}

\begin{document}

\def\ALERT#1{{\large\bf $\clubsuit$#1$\clubsuit$}}

\numberwithin{equation}{section}
\newtheorem{defin}{Definition}
\newtheorem{theorem}{Theorem}
\newtheorem{proposition}{Proposition}
\newtheorem{notice}{Notice}
\newtheorem{hypothesis}{Hypothesis}
\newtheorem{lemma}{Lemma}
\newtheorem{cor}{Corollary}
\newtheorem{example}{Example}
\newtheorem{remark}{Remark}
\newtheorem{conj}{Conjecture}
\def\begproof{\noindent{\bf Proof: }}
\def\endproof{\par\rightline{\vrule height5pt width5pt depth0pt}\medskip}
\def\div{\nabla\cdot}
\def\rot{\nabla\times}
\def\sign{{\rm sign}}
\def\arsinh{{\rm arsinh}}
\def\arcosh{{\rm arcosh}}
\def\diag{{\rm diag}}
\def\const{{\rm const}}
\def\d{\,\mathrm{d}}
\def\eps{\varepsilon}
\def\phi{\varphi}
\def\theta{\vartheta}
\def\N{\mathbb{N}}
\def\R{\mathbb{R}}
\def\C{\hbox{\rlap{\kern.24em\raise.1ex\hbox
      {\vrule height1.3ex width.9pt}}C}}
\def\P{\hbox{\rlap{I}\kern.16em P}}
\def\Q{\hbox{\rlap{\kern.24em\raise.1ex\hbox
      {\vrule height1.3ex width.9pt}}Q}}
\def\M{\hbox{\rlap{I}\kern.16em\rlap{I}M}}
\def\Z{\hbox{\rlap{Z}\kern.20em Z}}
\def\({\begin{eqnarray}}
\def\){\end{eqnarray}}
\def\[{\begin{eqnarray*}}
\def\]{\end{eqnarray*}}
\def\part#1#2{\frac{\partial #1}{\partial #2}}
\def\partk#1#2#3{\frac{\partial^#3 #1}{\partial #2^#3}} 
\def\mat#1{{D #1\over Dt}}
\def\dx{\nabla_x}
\def\dv{\nabla_v}
\def\grad{\nabla}
\def\Norm#1{\left\| #1 \right\|}
\def\pmb#1{\setbox0=\hbox{$#1$}
  \kern-.025em\copy0\kern-\wd0
  \kern-.05em\copy0\kern-\wd0
  \kern-.025em\raise.0433em\box0 }
\def\bar{\overline}
\def\lbar{\underline}
\def\fref#1{(\ref{#1})}
\def\half{\frac{1}{2}}
\def\oo#1{\frac{1}{#1}}

\def\tot#1#2{\frac{\d #1}{\d #2}} 
\def\laplace{\Delta}
\def\d{\,\mathrm{d}}
\def\N{\mathbb{N}}
\def\R{\mathbb{R}}
\def\supp{\mbox{supp }}
\def\eps{\varepsilon}
\def\phi{\varphi}

\def\E{\mathcal{E}}
\def\bbE{\mathbb{E}}

\def\Ikn{\I_{k_n}}
\def\sIkn{(\Ikn)_{n\in\N}}

\def\O{\mathcal{O}}
\def\A{\mathcal{A}}
\def\calG{\mathcal{G}}
\def\calF{\mathcal{F}}
\def\calX{\mathcal{X}}
\def\calY{\mathcal{Y}}
\def\calT{\mathcal{T}}
\def\calH{\mathcal{H}}
\def\calR{\mathcal{R}}
\def\bbX{\mathbb{X}}

\def\laplaceD{\laplace}

\def\comment#1{\textbf{[#1]}}


\centerline{{\huge Mathematical Analysis of a System}}
\centerline{{\huge for Biological Network Formation}}
\vskip 7mm


\vskip 5mm

\centerline{
{\large Jan Haskovec}\footnote{Mathematical and Computer Sciences and Engineering Division,
King Abdullah University of Science and Technology,
Thuwal 23955-6900, Kingdom of Saudi Arabia; 
{\it jan.haskovec@kaust.edu.sa}}\qquad
{\large Peter Markowich}\footnote{Mathematical and Computer Sciences and Engineering Division,
King Abdullah University of Science and Technology,
Thuwal 23955-6900, Kingdom of Saudi Arabia; 
{\it peter.markowich@kaust.edu.sa}}\qquad
{\large Benoit Perthame}\footnote{Sorbonne Universit\'es, UPMC Univ Paris 06, UMR 7598, Laboratoire Jacques-Louis Lions, F-75005, Paris, France, {\it benoit.perthame@upmc.fr}}${}^,$\footnote{CNRS, UMR 7598, Laboratoire Jacques-Louis Lions, F-75005, Paris, France}
}
\vskip 6mm


\noindent{\bf Abstract.}
Motivated by recent physics papers describing rules for natural network formation,
we study an elliptic-parabolic  system of partial differential equations proposed by Hu and Cai \cite{Hu, Hu-Cai}.
The model describes the pressure field thanks to Darcy's type equation and the dynamics of the conductance
network under pressure force effects with a diffusion rate $D$ representing randomness in the material structure. 

We prove the existence of global weak solutions and of local mild solutions and study their long term behavior.
It turns out that, by energy dissipation, steady states play a central role to understand the pattern capacity of the system.
We show that  for a large diffusion coefficient $D$, the zero steady state is stable.
Patterns occur for small values of $D$ because the zero steady state is Turing unstable in this range;
for $D=0$ we can exhibit a large class of dynamically stable (in the linearized sense) steady states. 
\vskip 5mm

\noindent{\bf Key words:} Pattern formation; Energy dissipation; Bifurcation analysis; Weak solutions; Stability; Turing instability; 
\\
\noindent{\bf Math. Class. No.:} 35K55; 35B32; 92C42

\section{Introduction}
\label{sec:introduction}

Network structures and dynamics, organization of leaf venation networks, vascular pattern formation
and their optimality in transport properties (electric, fluids, material) have been widely investigated
in the recent physics literature, in particular aiming at understanding natural networks.
Most of the methodological tools use discrete models  as in \cite{corson, couder, magnasco} to explain conductance dynamics.
Supply optimization has also been studied in more mathematical manner, see \cite{morel} and the references therein.
In opposition to the global effect of optimization, another explanation for networks topological structures,
purely local and based on mechanical laws has been proposed in \cite{Hu, Hu-Cai}.
Passing to the limit in the discrete model the authors also derive  a continuous model of network dynamics
which consists of a Poisson-type equation for the scalar pressure $p(t,x)$ coupled to a nonlinear diffusion equation
for the vector-valued conductance vector $m(t,x) $ of the network.  In \cite{Hu, Hu-Cai}, the authors propose the following system
with parameters $D\geq 0$, $c>0$,
\(
    -\div [(I + m\otimes m)\grad p] &=& S,   \label{eq1}\\
    \part{m}{t} - D^2\laplace m - c^2(m\cdot\grad p)\grad p + |m|^{2(\gamma-1)}m &=& 0, \label{eq2}
\)
on a bounded domain $\Omega\subset\R^d$, $d\leq 3$, with smooth boundary $\partial\Omega$,
subject to homogeneous Dirichlet boundary conditions on $\partial\Omega$ for $m$ and $p$,
\(   \label{BC_0}
   m(t,x) = 0,\qquad p(t,x)=0 \qquad \mbox{for } x\in\partial\Omega,\; t\geq 0.
\)
The system is completed with the initial condition for $m$,
\( \label{IC_0}
   m(t=0,x)=m^I(x)\qquad\mbox{for } x\in\Omega.
\)
The source term $S=S(x)$ is assumed to be independent of time and $\gamma \geq 1$ a parameter which is crucial for the type of networks formed \cite{Hu-Cai}.

The main mathematical interest of the PDE system for network formation stems from the highly unusual nonlocal coupling of the 
elliptic equation \eqref{eq1} for the pressure $p$ to the reaction-diffusion equation \eqref{eq2} for the conductance vector $m$
via the pumping term $+c^2 (\grad p\otimes\grad p) m$ and the latter term's potential equilibriation with the decay term
$- |m|^{2(\gamma-1)}m$. 
Also we remark that values of $\gamma<1$ (particularly $1/2 \leq \gamma <1$) make sense but in this paper we shall restrict to $\gamma\geq 1$.

Our purpose here is to study the existence of weak solution with $m(t) \in \left(H^1_0(\Omega)\right)^d$ (see Section\ref{sec:existence}).  The major difficulty being that a priori estimates  for $m$ are too weak to use elliptic regularizing effects for $p$.  Therefore weak solutions are just in the energy space that we recall below. However mild solutions with $m(t) \in \left(L^\infty(\Omega)\right)^d$ can be built (see Section~\ref{sec:mild_sols}) by a perturbation method. But, mostly, we wish to unravel some qualitative features which sustain the pattern formation properties of the system \eqref{eq1}--\eqref{eq2}. These rely on the existence of a large family of singular non-zero stationary solutions for $D=0$ which we build in Section~\ref{sec:pattern formation}.
\\

A major observation concerning system  \eqref{eq1}--\eqref{eq2} is  the energy-type functional
\(   \label{energy}
   \E(m) := \frac12 \int_\Omega \left(D^2 |\grad m|^2 + \frac{|m|^{2\gamma}}{\gamma} + c^2|m\cdot\grad p[m]|^2  + c^2|\grad p[m]|^2 \right) \d x,
\)
where $p=p[m] \in H^1_0(\Omega)$ is the unique solution of the Poisson equation \eqref{eq1} with $m$ given. We prove at first:

\begin{lemma}\label{lem:energy}
The energy \eqref{energy} is nonincreasing along smooth solutions of \eqref{eq1}--\eqref{eq2} and satisfy
\[
    \tot{}{t} \E(m) = - \int_\Omega \left( \part{m}{t}(t,x)\right)^2 \d x.
\] 
\end{lemma}

For weak solutions built in Theorem~\ref{thm:global_ex}, we recover as usual the weaker form
\[
 \E(m(t)) \leq  \E(m^I).
\]

\begproof
Multiplication of \eqref{eq1} by $p$ and integration by parts yields
\[
    \int_\Omega \left( |\grad p|^2 + |m\cdot\grad p|^2 - pS \right) \d x = 0.
\]
Subtracting the $c^2$-multiple of the above identity from \eqref{energy}, we obtain
\[
   \E(m(t)) = \frac12 \int_\Omega \left(D^2 |\grad m|^2 + \frac{|m|^{2\gamma}}{\gamma} - c^2|m\cdot\grad p|^2 - c^2|\grad p|^2 + 2c^2pS\right) \d x,
\]
so that, after integration by parts in suitable terms,
\[
   \tot{\E}{t}(m(t)) &=& - \int_\Omega D^2\laplace m \cdot \partial_t m \d x + \int_\Omega |m|^{2(\gamma-1)}m\cdot\partial_t m \d x
      - c^2 \int_\Omega (m\cdot\grad p) \grad p \cdot \partial_t m \d x \\
      && + c^2 \int_\Omega \grad\cdot \left[ (m\cdot\grad p) m\right]\partial_t p \d x + c^2\int_\Omega (\laplace p)(\partial_t p) \d x + c^2 \int_\Omega (\partial_t p) S \d x \\
      &=& - \int_\Omega \left[ D^2\laplace m - |m|^{2(\gamma-1)}m + c^2 (m\cdot\grad p) \grad p \right] \cdot \partial_t m \d x \\
      && + c^2 \int_\Omega \left[ \grad\cdot\bigl(\grad p + (m\otimes m)\grad p\bigr) + S\right] \partial_t p \d x \\
      &=& - \int_\Omega |\partial_t m|^2 \d x.
\]
\endproof

In the following, generic, not necessarily equal, constants will be denoted by $C$, but we will make specific use of the Poincar\'e constant  $C_\Omega$, i.e.,
\[
    \Norm{u}_{L^2(\Omega)} \leq C_\Omega \Norm{\grad u}_{L^2(\Omega)} \qquad\mbox{for all } u\in H_0^1(\Omega).
\]

As far as pattern formation properties of this system are concerned, let us point out that the sign of the term $-c^2(m\cdot\grad p)\grad p$ in \eqref{eq2} is crucial for structured networks  to be formed. Indeed, we have
\begin{lemma}\label{rem:non-pattern-formation}
For arbitrary functions $p$, the only weak stationary solution of the system \eqref{eq2}--\eqref{BC_0}
with $c^2$ instead of $-c^2$, i.e.,
\[
  - D^2\laplace m + c^2(m\cdot\grad p)\grad p + |m|^{2(\gamma-1)}m = 0\quad\mbox{in }\Omega,\qquad m=0\quad\mbox{on }\partial\Omega,
\]
is $m_0\equiv 0$.
\end{lemma}
\begproof
This follows directly from multiplication of the above equation
by $m$ and integration by parts,
\[
   D^2 \int_\Omega |\grad m|^2 \d x + c^2 \int_\Omega |m\cdot\grad p|^2 \d x + \int_\Omega |m|^{2\gamma} \d x = 0.
\]
Moreover, all weak solutions of equation \eqref{eq2} with $c^2$ instead of $-c^2$, i.e.,
\[
    \part{m}{t} - D^2\laplace m + c^2(m\cdot\grad p)\grad p + |m|^{2(\gamma-1)}m =0
\]
subject to homogeneous Dirichlet boundary condition and with arbitrary $\grad p$,
converge exponentially fast to zero in the $L^2$ sense.
Again, multiplication of the $m$-equation and integration by parts yields
\[
    \frac12 \tot{}{t} \int_\Omega |m|^2 \d x + D^2 \int_\Omega |\grad m|^2 \d x = -c^2 \int_\Omega |m\cdot\grad p|^2 \d x - \int_\Omega |m|^{2\gamma} \d x \leq 0.
\]
Therefore, using the Poincar\'e inequality with constant $C_\Omega$ gives
\[
    \frac12 \tot{}{t} \int_\Omega |m|^2 \d x \leq - D^2 \int_\Omega |\grad m|^2 \d x \leq - \frac{D^2}{C_\Omega^2} \int_\Omega |m|^2 \d x,
\]
and the Gronwall lemma provides the exponential convergence of the $L^2$-norm of $m$ to zero.
\endproof

\section{Global existence of weak solutions}
\label{sec:existence}

Our first goal is to prove the existence of global weak solutions of the system \eqref{eq1}--\eqref{IC_0}.
The main result is the following.

\begin{theorem}[Weak solutions]
Let $S\in L^2(\Omega)$ and $m^I \in H_0^1(\Omega)^d \cap L^{2\gamma}(\Omega)^d$.
Then the problem \eqref{eq1}--\eqref{IC_0} admits a global weak solution $(m,p[m])$ with $\E(m)\in L^\infty(0,\infty)$ , i.e., with
\[
     m \in L^\infty(0,\infty; H_0^1(\Omega)) \cap L^\infty(0,\infty; L^{2\gamma}(\Omega)), &\qquad&
    \partial_t m \in L^2((0,\infty)\times\Omega),\\
    \grad p \in L^\infty(0,\infty; L^2(\Omega)), &\qquad&
    m\cdot\grad p \in L^\infty(0,\infty; L^2(\Omega)).
\]
These solution satisfy the energy inequality,  with $\E$ given by \eqref{energy}, 
\(    \label{energy_ineq}
    \E(m(t)) + \int_0^t \int_\Omega \left( \part{m}{t}(s,x)\right)^2 \d x\d s \leq \E(m^I) \qquad\mbox{for all } t \geq 0.
\)
\label{thm:global_ex}
\end{theorem}

We proceed by proving existence of solutions of a regularized problem and a subsequent limit passage.
For this, we need some analytical results on an auxiliary parabolic and an elliptic problem,
which is the subject to the following two subsections.

\subsection{Analysis of auxiliary problems}
\label{subsec:auxiliary}

We first consider the semilinear parabolic problem on $\Omega$
\(   \label{m-aux}
    \part{m}{t} - D^2\laplace m + |m|^{2(\gamma-1)}m = f
\)
subject to the initial and boundary conditions
\(   \label{m-aux-IC-BC}
   m(t=0)=m^I\quad\mbox{in }\Omega,\qquad m=0\quad\mbox{on }\partial\Omega.
\)

\begin{lemma}\label{lem:m-aux}
For every $D>0$, $\gamma\geq 1$, $T>0$ and $f\in L^2((0,T)\times\Omega)^d$, the PDE \eqref{m-aux}--\eqref{m-aux-IC-BC}
with $m^I\in L^2(\Omega)^d$ admits a unique weak solution
$m\in L^2(0,T; H^1_0(\Omega))^d \cap C([0,T]; L^2(\Omega))^d \cap L^{2\gamma}((0,T)\times\Omega)^d$ with
\(  \label{m-aux-est1}
   \Norm{m}_{L^\infty(0,T; L^2(\Omega))} \leq \Norm{m^I}_{L^2(\Omega)} + C \Norm{f}_{L^2((0,T)\times\Omega)}
\)
and
\(  \label{m-aux-est2}
   \Norm{m}_{L^{2\gamma}((0,T)\times\Omega)} + \Norm{\grad m}_{L^2((0,T)\times\Omega)} \leq C \left( \Norm{m^I}_{L^2(\Omega)} + \Norm{f}_{L^2((0,T)\times\Omega)}  \right).
\)
Moreover, if $m^I\in H_0^1(\Omega)^d$, then $m\in L^\infty(0,T; H_0^1(\Omega))^d \cap L^2(0,T; H^2(\Omega))^d \cap L^\infty(0,T; L^{2\gamma}(\Omega))^d$
with $\partial_t m \in L^2((0,T)\times\Omega)^d$
and the estimates hold
\(
   \Norm{m}_{L^\infty(0,T; H_0^1(\Omega))} &\leq& C \Norm{f}_{L^2((0,T)\times\Omega)} + \Norm{m^I}_{H_0^1(\Omega)}, \label{m-aux-est4a}\\
   \Norm{\laplace m}_{L^2((0,T)\times\Omega)} &\leq& C \left( \Norm{f}_{L^2((0,T)\times\Omega)} + \Norm{m^I}_{H^1_0(\Omega)} \right), \label{m-aux-est4b}\\
   \Norm{m}_{L^\infty(0,T; L^{2\gamma}(\Omega))} &\leq& C \left( \Norm{f}_{L^2((0,T)\times\Omega)} + \Norm{m^I}_{H^1_0(\Omega)} 
             + \Norm{m^I}_{L^{2\gamma}(\Omega)} \right), \label{m-aux-est4c} \\
   \Norm{\partial_t m}_{L^2((0,T)\times\Omega)} &\leq& C \left( \Norm{f}_{L^2((0,T)\times\Omega)} + \Norm{m^I}_{H^1_0(\Omega)}
       + \Norm{m^I}_{L^{2\gamma}(\Omega)}\right). \label{m-aux-est5}
\)
The constants $C$ in \eqref{m-aux-est1}--\eqref{m-aux-est5} are independent of $f$ and $m^I$.
\end{lemma}

\begproof
A solution of \eqref{m-aux} can be constructed as a solution of the differential inclusion
\[
    \part{m}{t} + \partial I[m] \ni f
\]
with the functional $I: L^2(\Omega) \to [0,+\infty]$ given by
\[
    I[m] := \frac12 \int_\Omega |\grad m|^2 \d x + \frac{1}{2\gamma} \int_\Omega |m|^{2\gamma}\d x,\qquad\mbox{if } m\in H^1_0(\Omega)\cap L^{2\gamma}(\Omega),
\]
and $I[m]:= +\infty$ otherwise.
It can be easily checked that $I$ is proper with dense domain, strictly convex and lower semicontinuous, which implies the existence
of a unique solution $m\in L^2(0,T; H^1_0(\Omega))^d \cap L^{2\gamma}((0,T)\times\Omega)^d$.

Now let $f$ and $m^I$ be smooth. Multiplication of \eqref{m-aux} by $m$ and integration by parts yields
\[
    \frac12 \tot{}{t} \int_\Omega |m|^2 \d x + D^2\int_\Omega |\grad m|^2 \d x + \int_\Omega |m|^{2\gamma} \d x =
      \int_\Omega mf \d x \leq \frac12 \int_\Omega |m|^2 \d x + \frac12 \int_\Omega |f|^2 \d x.
\]
Using the Poincar\'e inequality, we have
\[
    \frac12 \tot{}{t} \int_\Omega |m|^2 \d x + \frac{D^2}{2} \int_\Omega |\grad m|^2 \d x + \int_\Omega |m|^{2\gamma} \d x
      \leq \frac{D^2}{2C_\Omega^2} \int_\Omega |f|^2 \d x,
\]
and \eqref{m-aux-est1}, \eqref{m-aux-est2} follow by a standard density argument.

Using $\laplace m$ as a test function, we obtain after an integration by parts
\[
    \frac12 \tot{}{t} \int_\Omega |\grad m|^2 \d x + D^2\int_\Omega |\laplace m|^2 \d x
       - \int_\Omega (|m|^{2(\gamma-1)}m)\cdot\laplace m \d x = \int_\Omega f\laplace m \d x.
\]
We integrate by parts in the term
\[
   - \int_\Omega (|m|^{2(\gamma-1)}m)\cdot\laplace m \d x &=&
       \int_\Omega \sum_{i,j=1}^d \part{(|m|^{2(\gamma-1)}m_i)}{x_j} \part{m_i}{x_j} \d x \\
      &=& 2(\gamma-1) \int_\Omega |m|^{2(\gamma-2)} \sum_{i,j=1}^d \left(m_i \part{m_i}{x_j}\right)^2 \d x
         + \int_\Omega |m|^{2(\gamma-1)} |\grad m|^2 \d x \\ &\geq& 0.
\]
Therefore, we have
\[
    \frac12 \tot{}{t} \int_\Omega |\grad m|^2 \d x + D^2 \int_\Omega |\laplace m|^2 \d x &\leq&
          \int_\Omega f\laplace m \d x \\
      &\leq& \frac1{2D^2} \int_\Omega |f|^2 \d x + \frac{D^2}2 \int_\Omega |\laplace m|^2 \d x,
\]
which after integration in time and application of another density argument
gives the estimates \eqref{m-aux-est4a} and \eqref{m-aux-est4b}.

Next, we define the energy
\[
    E := \frac{D^2}2 \int_\Omega |\grad m|^2 \d x + \frac{1}{2\gamma} \int_\Omega |m|^{2\gamma} \d x
\]
and calculate, for smooth solutions,
\[
    \tot{E}{t} &=& - D^2 \int_\Omega \partial_t m \cdot \laplace m \d x + \int_\Omega |m|^{2(\gamma-1)}m\cdot\partial_t m \d x \\
      &=& \int_\Omega \partial_t m \cdot(-\partial_t m + f) \d x.
\]
This implies, after integration in time,
\[
    E(T) + \int_0^T \int_\Omega |\partial_t m|^2 \d x\d t &=& E(0) + \int_0^T \int_\Omega f\cdot\partial_t m \d x\d t \\
       &\leq& E(0) + \frac12 \int_0^T \int_\Omega |f|^2 \d x\d t + \frac12 \int_0^T \int_\Omega |\partial_t m|^2 \d x\d t,
\]
and, again by a density argument, \eqref{m-aux-est4c}, \eqref{m-aux-est5} follow.
\endproof

We will also need the following Lemma concerning the algebraic term $|m|^{2(\gamma-1)}m$.

\begin{lemma}\label{lem:m-alg-conv}
Fix $\gamma\geq 1$ and let the sequence $\{m^k\}_{k\in\N}$ be uniformly bounded in $L^{2\gamma}((0,T)\times\Omega)$
and converging to $m$ in the norm topology of $L^2((0,T)\times\Omega)$ as $k\to\infty$.
Then for every test function $\phi\in C^\infty_c([0,T)\times\Omega)$,
\[
   \int_0^T \int_\Omega |m^k|^{2(\gamma-1)}m^k \phi \d x\d t \quad\rightarrow\quad \int_0^T \int_\Omega |m|^{2(\gamma-1)}m\, \phi \d x\d t
     \qquad\mbox{as } k\to\infty.
\]
\end{lemma}

\begproof
We construct the Young measure $\nu_{x,t}$ on $\R^d$ corresponding to the weakly converging sequence $m^k$ in $L^{2\gamma}((0,T)\times\Omega)$,
so that for every function $h\in C_0(\R^d)$ and every test function $\phi\in C^\infty_c([0,T)\times\Omega)$,
\(    \label{YoungMConv}
    \int_0^T \int_\Omega h(m^k)\phi \d x\d t \quad\stackrel[k\to\infty]{}{\longrightarrow}\quad
       \int_0^T \int_\Omega \left(\int_{\R^d} h(y) \d\nu_{x,t}(y) \right) \phi(t,x) \d x\d t.
\)
However, due to the assumed strong convergence of $m^k$ in $L^2((0,T)\times\Omega)$, 
we have $h(m^k)\to h(m)$ strongly in $L^2((0,T)\times\Omega)$ if $h\in C_0(\R^d)$.
Consequently, $\nu_{x,t}(y) = \delta(y - m(t,x))$.
Then, due to the Fundamental Theorem of Young measures \cite{Fonseca-Leoni}, \eqref{YoungMConv}
holds also for $h\in C(\R^d)$ if the sequence $h(m^k)$ is equiintegrable.
In the particular case $h(s) = |s|^{2(\gamma-1)}s$ the equiintegrability is verified by
\[
    \int_{\{|h(m^k)|\geq K\}} |h(m^k)| \d x\d t =  \int_{\{|h(m^k)|\geq K\}} |m^k|^{2\gamma} |m^k|^{-1}\d x\d t 
      \leq \Norm{m^k}_{L^{2\gamma}((0,T)\times \Omega)}^{2\gamma} K^{-1/({2\gamma-1})} \leq C K^{-1/({2\gamma-1})},
\]
for every $K>0$.
\endproof

Next, we study the properties of solutions of the regularized Poisson equation
\(   \label{Poisson_reg}
   -\div [\grad p + \bar m (\bar m\cdot\grad p)\ast\eta ] = S,
\)
where $\bar m$ and $S$ are given functions on $\Omega$
and $\eta=\eta(|x|)\in C^\infty(\R^d)\cap L^\infty(\R^d)$ is a smooth mollifier
with real, nonnegative Fourier transform (in particular, we shall use the heat kernel later on).
The convolution $f\ast\eta$ for $f\in L^1(\Omega)$ is defined as
\[
   f\ast\eta(x) = \int_{\R^d} f(y)\eta(x-y) \d y,
\]
where we extend $f$ by zero to $\R^d$, i.e., $f(y)=0$ for $y\in\R^d\setminus\Omega$.
We need the following technical Lemma.

\begin{lemma}\label{lem:nonneg_convolution}
For any $u\in L^1(\R^d)$ and $\eta\in L^\infty(\R^d) \cap L^1(\R^d)$ with nonnegative Fourier transform $\hat\eta \geq 0$ on $\R^d$,
the identity holds
\[
    \int_{\R^d} (u\ast\eta)u \d x = \int_{\R^d} |u\ast\rho|^2 \d x \geq 0,
\]
where $\rho$ is the inverse Fourier transform of $(\hat\eta)^{1/2}$.
\end{lemma}

\begproof
With the Parseval's identity we have
\[
    \int_{\R^d} (u\ast\eta)u \d x = \int_{\R^d} \widehat{u\ast\eta}(\xi)\; \bar{\hat u(\xi)} \d\xi =
    \int_{\R^d} \hat u(\xi)\hat\eta(\xi) \bar{\hat u(\xi)} \d\xi
\]
and since, by assumption, $\hat\eta = |\hat\rho|^2 = \hat\rho\bar{\hat\rho}$, this is further equal to
\[
   \int_{\R^d} [\hat u(\xi)\hat\rho(\xi)] \bar{[\hat u(\xi)\hat\rho(\xi)]} \d\xi =
     \int_{\R^d} |u\ast\rho|^2 \d x \geq 0.
\]
\endproof

\begin{lemma}\label{lem:est_gradp}
For every $\bar m\in L^2(\Omega)$ and $S\in L^2(\Omega)$, the regularized Poisson equation \eqref{Poisson_reg}
has a unique weak solution $p\in H^1_0(\Omega)$.
Moreover,
\(  \label{est_gradp}
    \Norm{\grad p}_{L^2(\Omega)} \leq C_\Omega\Norm{S}_{L^2(\Omega)}.
\)
\end{lemma}

\begproof
We define the bilinear form $B: H^1_0(\Omega) \times H^1_0(\Omega) \to \R$,
\[
    B(p,\phi) := \int_\Omega \grad p\cdot\grad\phi \d x + \int_\Omega [(\bar m\cdot\grad p)\ast\eta](\bar m\cdot\grad\phi) \d x.
\]
Then, due to Lemma~\ref{lem:nonneg_convolution},
\[
    \int_\Omega [(\bar m\cdot\grad p)\ast\eta](\bar m\cdot\grad p) \d x \geq 0,
\]
so that we have the ellipticity condition $B(p,p) \geq \Norm{\grad p}^2_{L^2(\Omega)}$.
Continuity of $B$ follows from
\[
   \int_\Omega [(\bar m\cdot\grad p)\ast\eta](\bar m\cdot\grad\phi) \d x &\leq&
      \Norm{(\bar m\cdot\grad p)\ast\eta}_{L^\infty(\Omega)} \Norm{\bar m\cdot\grad\phi}_{L^1(\Omega)} \\
   &\leq& \Norm{\eta}_{L^\infty(\Omega)} \Norm{\bar m\cdot\grad p}_{L^1(\Omega)} \Norm{\bar m\cdot\grad\phi}_{L^1(\Omega)} \\
   &\leq& C \Norm{\bar m}_{L^2(\Omega)}^2 \Norm{\grad p}_{L^2(\Omega)} \Norm{\grad\phi}_{L^2(\Omega)}.
\]
Therefore, the Lax-Milgram theorem provides the existence of unique weak solutions $p\in H^1_0(\Omega)$ of \eqref{Poisson_reg}.
Then, \eqref{est_gradp} follows from
\[
   \Norm{\grad p}^2_{L^2(\Omega)} \leq B(p,p) = \int_\Omega S p \d x \leq C_\Omega \Norm{S}_{L^2(\Omega)} \Norm{\grad p}_{L^2(\Omega)}.
\]
\endproof

\begin{lemma}\label{lem:weak-strong}
Let $\bar m^k$ be a sequence of functions converging to $\bar m$ in the norm topology of $L^2(\Omega)$, and denote by $p^k\in H^1_0(\Omega)$
the corresponding weak solutions of \eqref{Poisson_reg} subject to homogeneous Dirichlet boundary conditions.
Then $p^k$ converges strongly in $H^1(\Omega)$ to the unique solution $p$ of \eqref{Poisson_reg} as $k\to\infty$.
\end{lemma}

\begproof
Lemma~\ref{lem:est_gradp} provides the uniform bound $\Norm{\grad p^k}_{L^2(\Omega)} \leq C$, so that,
for a subsequence again denoted by $\grad p^k$, we have weak convergence $\grad p^k \rightharpoonup \grad p$ in $L^2(\Omega)$
for some $p\in H^1_0$. Moreover, denoting $q^k := (\bar m^k\cdot\grad p^k)\ast\eta$, we have
\[
   \Norm{q^k}_{L^\infty(\Omega)} \leq \Norm{\eta}_{L^\infty(\R^d)}\Norm{\bar m^k\cdot\grad p^k}_{L^1(\Omega)}
     \leq \Norm{\eta}_{L^\infty(\R^d)} \Norm{\bar m^k}_{L^2(\Omega)} \Norm{\grad p^k}_{L^2(\Omega)} \leq C,
\]
so that a subsequence of $q^k$ converges to $q$ weakly* in $L^\infty(\Omega)$ as $k\to\infty$.
Moreover, due to the strong convergence of $\bar m^k$ to $\bar m$ and weak convergence
of $\grad p^k$ to $\grad p$ in $L^2(\Omega)$, we have for every test function $\phi\in C^\infty_0(\Omega)$
\[
    \int_\Omega q^k\phi \d x &=& \int_{\R^d} (\bar m^k\cdot\grad p^k)(\phi\ast\eta) \d x  \\
      &\stackrel[k\to\infty]{}{\longrightarrow}& \int_{\R^d} (\bar m\cdot\grad p)(\phi\ast\eta) \d x 
      = \int_\Omega [(\bar m\cdot\grad p)\ast\eta] \phi \d x,
\]
and the limit $q$ is identified as $(\bar m\cdot\grad p)\ast\eta$.
Consequently, we can pass to the limit in the weak formulation of \eqref{Poisson_reg},
and using $p$ as a test function then yields
\[
    \int_\Omega |\grad p|^2 + (\bar m\cdot\grad p) q \d x = \int_\Omega S p \d x.     
\]
On the other hand, using $p^k$ as a test function in \eqref{Poisson_reg} for $p^k$, we have
\(    \label{tricky_eq}
     \lim_{k\to\infty} \int_\Omega |\grad p^k|^2 +  (\bar m^k\cdot\grad p^k) q^k \d x  &=& \lim_{k\to\infty} \int_\Omega S p^k \d x \\
       &=& \int_\Omega |\grad p|^2 + (\bar m\cdot\grad p) q \d x. \nonumber
\)
Now, we write
 \[
     \liminf_{k\to \infty} \int_\Omega |\grad p^k|^2 &\leq& 
     \limsup_{k\to \infty} \int_\Omega |\grad p^k|^2  \\ 
     &=& 
     \limsup_{k\to \infty} \int_\Omega \left[ |\grad p^k|^2 + (\bar m^k\cdot\grad p^k) q^k
           - (\bar m^k\cdot\grad p^k) q^k \right] \d x \\
     &\leq& 
     \limsup_{k\to \infty} \int_\Omega |\grad p^k|^2 + (\bar m^k\cdot\grad p^k) q^k \d x
         + \limsup_{k\to \infty} \int_\Omega - (\bar m^k\cdot\grad p^k) q^k \d x.
\]
Due to \eqref{tricky_eq}, the last line is equal to
\[
     \int_\Omega |\grad p|^2 + (\bar m\cdot\grad p) q \d x
        - \liminf_{k\to \infty} \int_\Omega (\bar m^k\cdot\grad p^k) q^k \d x.  
\]
Following Lemma~\ref{lem:nonneg_convolution}, we write the second integral as
\[
    \int_\Omega (\bar m^k\cdot\grad p^k) q^k \d x = \int_\Omega (\bar m^k\cdot\grad p^k) [(\bar m^k\cdot\grad p^k)\ast\eta] \d x 
       = \int_{\R^d} \left|(\bar m^k\cdot\grad p^k)\ast\rho\right|^2 \d x,
\]
with $\widehat{\eta} = |\widehat{\rho}|^2$. Clearly, the sequence $(\bar m^k\cdot\grad p^k)\ast\rho$ is bounded
in $L^2(\R^d)$ and a subsequence converges weakly to $(\bar m\cdot\grad p)\ast\rho$. 
Then the weak lower semicontinuity of the $L^2$ norm implies
\[
     - \liminf_{k\to \infty} \int_{\R^d} \left|(\bar m^k\cdot\grad p^k)\ast\rho\right|^2 \d x
        &\leq& - \int_{\R^d} \left|(\bar m\cdot\grad p)\ast\rho\right|^2 \d x \\
        &=& -\int_\Omega (\bar m\cdot\grad p) q \d x,
\]
so that, finally,
\[
     \liminf_{k\to \infty} \int_\Omega |\grad p^k|^2 \leq
     \int_\Omega |\grad p|^2 + (\bar m\cdot\grad p) q \d x - \int_\Omega (\bar m\cdot\grad p) q \d x
     = \int_\Omega |\grad p|^2 \d x.
\]
This directly implies $\lim_{k\to\infty} \Norm{\grad p^k}_{L^2(\Omega)} = \Norm{\grad p}_{L^2(\Omega)}$,
and thus the convergence of $\grad p^k$ is strong in $L^2(\Omega)$.
\endproof

\subsection{Existence of global solutions of the regularized problem}
\label{subsec:ex_reg}

For $\eps>0$ small we consider the perturbed problem
\(
    -\div [\grad p + m(m\cdot\grad p)\ast\eta_\eps ] &=& S,   \label{eq_pert1}\\
    \part{m}{t} - D^2\laplace m - c^2[(m\cdot\grad p)\ast\eta_\eps]\grad p + |m|^{2(\gamma-1)}m &=& 0, \label{eq_pert2}
\)
with $(\eta_\eps)_{\eps>0}$ the $d$-dimensional heat kernel $\eta_\eps(x) = (4\pi\eps)^{-d/2} \exp(-|x|^2/4\eps)$,
on a bounded domain $\Omega\subset\R^d$, $d\leq 3$, with smooth boundary $\partial\Omega$, subject to homogeneous Dirichlet boundary conditions
on $\partial\Omega$ for $m$ and $p$,
\(   \label{BC_pert1}
   m(t,x) = 0,\qquad p(t,x)=0 \qquad \mbox{for all } x\in\partial\Omega,\; t\geq 0,
\)
and the initial condition for $m$,
\(   \label{IC_pert1}
   m(t=0,x)=m^I(x) \qquad \mbox{for all } x\in\Omega.
\)
Let us note that the heat kernel $\eta_\eps$ satisfies the assumptions of Lemma~\ref{lem:nonneg_convolution} for any $\eps>0$.

\begin{theorem} [Existence for the perturbed model]
Let $m^I\in L^2(\Omega)^d$.
For any $\eps>0$ there exists a weak solution $(m,p)$ of the system \eqref{eq_pert1}--\eqref{IC_pert1}
with $p\in L^\infty(0,T; H^1_0(\Omega))$ and $m\in L^2(0,T; H^1_0(\Omega))^d$,
$\laplace m\in L^2((0,T)\times\Omega)^d$ and $\partial_t m \in L^2((0,T)\times\Omega)^d$.
\label{thm:ex_pert}
\end{theorem}

\begproof
We shall employ the Leray-Schauder fixed point theorem. We fix $\eps>0$ and construct the mapping $\Phi: \bar m\mapsto m$ in two steps:
For a given $\bar m \in L^2((0,T)\times\Omega)$ we set $p$
to be the unique weak solution of
\(   \label{Poisson_reg2}
   -\div [ \grad p + \bar m(\bar m\cdot\grad p)\ast\eta_\eps ] = S,
\)
constructed in Lemma~\ref{lem:est_gradp} (to be precise, we use a straightforward modification of Lemma~\ref{lem:est_gradp}
where the Lax-Milgram Theorem is applied in the space $L^2(0,T; H^1_0(\Omega))$).
By the same Lemma, we have the a priori estimate
$\Norm{\grad p}_{L^\infty(0,T; L^2(\Omega))} \leq C$ independent of $\bar m$.
We set $q_\eps := (\grad p\cdot\bar m)\ast \eta_\eps$ and note that
it is a priori bounded in $L^2(0,T; L^\infty(\Omega))$ due to
\(   \label{est_q}
    \Norm{q_\eps}_{L^2(0,T; L^\infty(\Omega))} \leq \Norm{\eta_\eps}_{L^\infty(\R^d)}\Norm{\grad p\cdot\bar m}_{L^2(0,T; L^1(\Omega))}
       \leq C_\eps \Norm{\grad p}_{L^\infty(0,T; L^2(\Omega))} \Norm{\bar m}_{L^2((0,T)\times\Omega)}.
\)
Then we employ Lemma~\ref{lem:m-aux} with $f:=c^2 q_\eps\grad p \in L^2((0,T)\times\Omega)$ and
set $\Phi(\bar m):=m \in L^2(0,T; H_0^1(\Omega))^d$ to be the unique weak solution of
\(   \label{m-pert}
   \part{m}{t} - D^2\laplace m + |m|^{2(\gamma-1)}m = c^2 q_\eps\grad p \label{eq_pert4}
\)
subject to the initial condition $m(t=0)=m^I$.
The Lemma provides the estimate
\[
    \Norm{m}_{L^\infty(0,T; L^2(\Omega))} + \Norm{\grad m}_{L^2((0,T)\times\Omega)} &\leq&
       C \left( \Norm{m^I}_{L^2(\Omega)} + c^2 \Norm{q_\eps\grad p}_{L^2((0,T)\times\Omega)} \right) \\
      &\leq& C \left( \Norm{m^I}_{L^2(\Omega)} + c^2 \Norm{q_\eps}_{L^2(0,T; L^\infty(\Omega))} \Norm{\grad p}_{L^\infty(0,T;L^2(\Omega))} \right) \\
      &\leq& \widetilde C + C_\eps \Norm{\grad p}^2_{L^\infty(0,T;L^2(\Omega))} \Norm{\bar m}_{L^2((0,T)\times\Omega)}
\]
for suitable constants $\widetilde C$, $C_\eps>0$.
This also implies an a priori bound on $\partial_t m$ in $L^2(0,T;H^{-1}(\Omega))$,
so that $\Phi: \bar m \mapsto m$ maps bounded sets in $L^2((0,T)\times\Omega)$
onto relatively compact ones.

To prove continuity of $\Phi$, consider a sequence $\bar m^k$ converging to $\bar m$
in the norm topology of $L^2((0,T)\times\Omega)$.
Due to Lemma~\ref{lem:est_gradp}, the sequence $\grad p^k$ of the corresponding solutions of \eqref{Poisson_reg2}
converges weakly-* in $L^\infty(0,T;L^2(\Omega))$ to some $\grad p$.
The bound \eqref{est_q} allows to extract a subsequence of $q_\eps^k := (\bar m^k\cdot\grad p^k)\ast\eta_\eps$
converging weakly in $L^2(0,T; L^\infty(\Omega))$ to $q_\eps := (\bar m\cdot\grad p)\ast\eta_\eps$.
Therefore, we can pass to the limit in the weak formulation of the regularized Poisson equation \eqref{Poisson_reg2}
and conclude that $p$ is its unique solution corresponding to $\bar m$.
To pass to the limit in \eqref{m-pert}, we use strong convergence of $\grad p^k$
in $L^2((0,T)\times\Omega)$ provided by Lemma~\ref{lem:weak-strong}
(strictly speaking, by its straightforward modification for time dependent functions $\bar m^k$).
Then the limit passage in the term $q_\eps^k\grad p^k$ is straightforward.
For the algebraic term $|m^k|^{2(\gamma-1)}m^k$ we employ Lemma~\ref{lem:m-alg-conv}
and finally conclude, by the uniqueness of solutions of \eqref{m-pert}, the continuity $m=\Phi(\bar m)$.

Finally, we have to prove that the set $\calY := \{m\in L^2((0,T)\times\Omega); \kappa\Phi(m) = m \mbox{ for some } 0 < \kappa\leq 1\}$ is bounded.
Note that the elements of this set are weak solutions of the system
\[
    -\div [ \grad p + m(m\cdot\grad p)\ast\eta_\eps ] &=& S,   \\
    \part{m}{t} - D^2\laplace m - \kappa c^2[(m\cdot\grad p)\ast\eta_\eps]\grad p + \kappa^{-2(\gamma-1)} |m|^{2(\gamma-1)}m &=& 0,
\]
subject to homogeneous Dirichlet boundary conditions for $m$ and $p$ on $\partial\Omega$ and the initial condition $m(t=0,x)=\kappa m^I(x)$ in $\Omega$.
Multiplication of the first equation by ${c^2}{\kappa}p$ and of the second equation by $m$, integration by parts
and subtraction of the two identities yields
\[
   \tot{}{t} \int_\Omega |m|^2 \d x + D^2\int_\Omega |\grad m|^2 \d x + \kappa^{-2(\gamma-1)} \int_\Omega |m|^{2\gamma}\d x =
   {c^2}{\kappa} \left( \int_\Omega pS \d x - \int_\Omega |\grad p|^2 \d x \right).
\]
This implies, for any $0 < \kappa\leq 1$,
\[
   \tot{}{t} \int_\Omega |m|^2 \d x + D^2\int_\Omega |\grad m|^2 \d x
   &\leq& {c^2}{\kappa} \left( C_p \Norm{S}_{L^2(\Omega)} \Norm{\grad p}_{L^2(\Omega)} - \Norm{\grad p}_{L^2(\Omega)}^2 \right) \\
   &\leq& {c^2} C_\Omega^2 \Norm{S}^2_{L^2(\Omega)},
\]
where we used the Poincar\'e inequality with the constant $C_\Omega$.
This immediately gives the a priori boundedness of $m$ in $L^\infty(0,T; L^2(\Omega))$
and thus the boundedness of the set $\calY$.

Finally, estimates \eqref{m-aux-est4b} and \eqref{m-aux-est5} of Lemma~\ref{lem:m-aux}
with $f:= q_\eps\grad p \in L^2((0,T)\times\Omega)^d$ imply that
$\laplace m\in L^2((0,T)\times\Omega)^d$ and $\partial_t m \in L^2((0,T)\times\Omega)^d$.
\endproof

\subsection{The limit $\eps\to 0$ in System \eqref{eq_pert1}--\eqref{eq_pert2}}

We shall now pass to the limit $\eps\to 0$ in \eqref{eq_pert1}--\eqref{IC_pert1} and obtain a global solution
of the system \eqref{eq1}--\eqref{IC_0}. The main tool is the dissipation of the modified energy
\(   \label{energy_eps}
   \E_\eps(m) := \frac12 \int \left(D^2 |\grad m|^2 + \frac{|m|^{2\gamma}}{\gamma} + c^2(m\cdot\grad p)[(m\cdot\grad p)\ast\eta_\eps]  + c^2|\grad p|^2 \right) \d x.
\)
Note that by Lemma~\ref{lem:nonneg_convolution} we have
\[
   \int_\Omega (m\cdot\grad p)[(m\cdot\grad p)\ast\eta_\eps] \d x = \int_{\R^d} |(m\cdot\grad p)\ast\rho_\eps|^2 \d x \geq 0
\]
with $\hat \eta_\eps = |\hat \rho_\eps|^2$, so that $\E_\eps(m(t)) \geq 0$.

\begin{lemma}\label{lem:energy_eps}
Let $(m,p)$ be a solution of \eqref{eq_pert1}--\eqref{eq_pert2} constructed in Theorem~\ref{thm:ex_pert}
and assume $\E_\eps(m^I) < \infty$.
Then the energy \eqref{energy_eps} satisfies
\(    \label{energy_eps_ineq}
    \E_\eps(m(t)) + \int_0^t \int_\Omega \left( \part{m}{t}(s,x)\right)^2 \d x\d s = \E_\eps(m^I) \qquad\mbox{for all } 0 \leq t \leq T.
\)
\end{lemma}

\begproof
Denote $q_\eps:=(m\cdot\grad p)\ast\eta_\eps$. Recalling that $q_\eps\in L^2(0,T;L^\infty(\Omega))$
and $m\cdot\grad p \in L^2(0,T;L^1(\Omega))$, we multiply \eqref{eq_pert1} by $p$ and integrate by parts, obtaining
\[
    \int_\Omega |\grad p|^2 + (m\cdot\grad p)q_\eps \d x = \int_\Omega pS  \d x.
\]
Subtracting the $c^2$-multiple of the above identity from \eqref{energy_eps}, we obtain
\[
   \E_\eps(t) = \frac12 \int \left(D^2 |\grad m|^2 + \frac{|m|^{2\gamma}}{\gamma} - c^2(m\cdot\grad p)q_\eps - c^2|\grad p|^2 + 2c^2pS\right) \d x,
\]
so that, after integration by parts in the suitable terms (note that Theorem~\ref{thm:ex_pert} provides
enough regularity on $m$ and $p$ for the below calculation), we have
\[
   \tot{\E_\eps(t)}{t} &=& - \int D^2\laplace m \cdot \partial_t m \d x + \int |m|^{2(\gamma-1)}m\cdot\partial_t m \d x
      - c^2 \int q_\eps \grad p \cdot \partial_t m \d x \\
      && + c^2 \int \grad\cdot \left( q_\eps m\right)\partial_t p \d x + c^2\int (\laplace p)(\partial_t p) \d x + c^2 \int (\partial_t p) S \d x \\
      &=& - \int \left[ D^2\laplace m - |m|^{2(\gamma-1)}m + c^2 q_\eps \grad p \right] \cdot \partial_t m \d x \\
      && + c^2 \int \left[ \grad\cdot\bigl(\grad p + q_\eps m\bigr) + S\right] \partial_t p \d x \\
      &=& - \int |\partial_t m|^2 \d x,
\]
where we used the identity
\[
   \int_\Omega \partial_t [(m\cdot\grad p)q_\eps] \d x = 2 \int_\Omega q_\eps\grad p\cdot\partial_t m \d x
      + 2 \int_\Omega q_\eps m\cdot\grad (\partial_t p) \d x,
\]
which holds due to the symmetry of the kernel $\eta_\eps(x) = \eta_\eps(|x|)$.
Integration of the above result in time yields \eqref{energy_eps_ineq}.
\endproof

We are now ready to pass to the limit $\eps\to 0$ in \eqref{eq_pert1}--\eqref{eq_pert2}.

\begin{lemma}\label{lem:eps_to_zero}
Let $(m^\eps,p^\eps)_{\eps>0}$ be a family of weak solution of \eqref{eq_pert1}--\eqref{eq_pert2}
constructed in Theorem~\ref{thm:ex_pert} and assume $\E(m^I) < \infty$. Then there exists a subsequence converging
to $(m,p)$ as $\eps\to 0$, where $(m,p)$ is a weak solution of \eqref{eq1}--\eqref{IC_0}, satisfying
the energy dissipation inequality \eqref{energy_ineq}.
\end{lemma}

\begproof
Note that the Poisson equation \eqref{eq_pert1} at $t=0$ implies
\(    \label{reg_Poisson_identity}
    \int_\Omega |\grad p^\eps[m^I]|^2 \d x + \int_\Omega |(m^I\cdot\grad p^\eps[m^I])\ast\rho_\eps|^2 \d x = \int_\Omega p^\eps[m^I] S \d x,
\)
such that $\E_\eps(m^I)$ is uniformly bounded as $\eps\to 0$.
Then the energy dissipation given by Lemma~\ref{lem:energy_eps} provides the following uniform a priori estimates,
\[
    m^\eps \in L^\infty(0,T; H_0^1(\Omega)) \cap L^\infty(0,T; L^{2\gamma}(\Omega)), &\qquad&
    \partial_t m^\eps \in L^2((0,T)\times\Omega),\\
    \grad p^\eps \in L^\infty(0,T; L^2(\Omega)),&\qquad&
    (m^\eps\cdot\grad p^\eps)\ast\rho_\eps \in L^\infty(0,T; L^2(\R^d)),
\]
with $\hat \eta_\eps = |\hat \rho_\eps|^2$. The last bound implies
a uniform estimate on $q_\eps := (m^\eps\cdot\grad p^\eps)\ast\eta_\eps$ in $L^\infty(0,T; L^2(\Omega))$.
Indeed, taking any test function $\phi\in L^1(0,T; L^2(\Omega))$, we have
\[
    \int_0^T \int_\Omega q_\eps \phi \d x\d t &=& \int_0^T \int_{\R^d} [(m^\eps\cdot\grad p^\eps)\ast\rho_\eps][\phi\ast\rho_\eps] \d x\d t \\
     &\leq& \int_0^T \Norm{(m^\eps\cdot\grad p^\eps)\ast\rho_\eps}_{L^2(\R^d)} \Norm{\phi\ast\rho_\eps}_{L^2(\R^d)} \d t \\
     &\leq& \Norm{(m^\eps\cdot\grad p^\eps)\ast\rho_\eps}_{L^\infty(0,T; L^2(\R^d))} \int_0^T \Norm{\phi}_{L^2(\Omega)}\Norm{\rho_\eps}_{L^1(\R^d)} \d t \\
     &=& \Norm{(m^\eps\cdot\grad p^\eps)\ast\rho_\eps}_{L^\infty(0,T; L^2(\R^d))} \Norm{\phi}_{L^1(0,T; L^2(\Omega))},     
\]
where we used the fact that, by definition of $\rho_\eps$, $\Norm{\rho_\eps}_{L^1(\R^d)}=1$ for every $\eps>0$.
Therefore, by duality, $q_\eps$ is uniformly bounded in $L^\infty(0,T; L^2(\Omega))$ and there exists
a subsequence converging weakly-* in this space to some $q\in L^\infty(0,T; L^2(\Omega))$.
We note that due to the compact embedding (Corollary 4 in \cite{Simon87}), a subsequence of $m^\eps$ converges to some $m$
in the norm topology of $L^2((0,T)\times\Omega)$. Then, a slight modification of Lemma~\ref{lem:weak-strong}
provides the strong convergence of $p^\eps$ to $p$ in $L^2(0,T; H_0^1(\Omega))$,
where $p$ is the unique solution of the Poisson equation \eqref{eq1} with $m$.
Consequently, the product $m^\eps\cdot\grad p^\eps$ converges strongly to $m\cdot\grad p$ in $L^1((0,T)\times\Omega)$,
and for every test function $\phi\in C^\infty_0([0,T)\times\Omega)$ we have
\[
    \int_0^T \int_{\R^d} q_\eps\phi \d x\d t = \int_0^T \int_{\R^d} (m^\eps\cdot\grad p^\eps) (\phi\ast\eta_\eps) \d x
      \;\stackrel[\eps\to 0]{}{\longrightarrow}\; \int_0^T \int_{\R^d} (m\cdot\grad p)\, \phi \d x,
\]
where we used the fact that $\phi\ast\eta_\eps$ converges to $\phi$ in $C([0,T]\times\bar\Omega)$ as $\eps\to 0$
due to the Arzela-Ascoli theorem.
Therefore, we identify the limit $q=m\cdot\grad p$.

We are now ready to pass to the limit in the weak formulation of the nonlinear terms of \eqref{eq_pert1}--\eqref{eq_pert2}.
The term $q_\eps m^\eps$ in \eqref{eq_pert1} converges to $(m\cdot\grad p)m$ due to the weak-* convergence
of $q_\eps$ in $L^\infty(0,T; L^2(\Omega))$ and strong convergence of $m^\eps$ in $L^2((0,T)\times\Omega)$.
The term $q_\eps\grad p^\eps$ in \eqref{eq_pert2} converges to $(m\cdot\grad p)\grad p$ due to the strong convergence
of $\grad p^\eps$ in $L^2((0,T)\times\Omega)$.
Finally, the limit passage in the term $|m^\eps|^{2(\gamma-1)}m^\eps$ is provided by Lemma~\ref{lem:m-alg-conv}
due to the uniform boundedness of $m^\eps$ in $L^{2\gamma}((0,T)\times\Omega)$.

The energy dissipation inequality \eqref{energy_ineq} follows from \eqref{energy_eps_ineq}
due to the weak lower semicontinuity of terms defining $\E(m)$
and from the fact that $\E_\eps(m^I)\to \E(m^I)$ as $\eps\to 0$.
Indeed, Lemma \eqref{lem:weak-strong} provides strong converge of $\grad p^\eps[m^I]$ to $\grad p[m^I]$
in $L^2(\Omega)^d$. Due to the embedding of $H_0^1(\Omega)$ into $L^6(\Omega)$ for $d\leq 3$,
the term $m^I\cdot\grad p^\eps[m^I]$ converges strongly in $L^{3/2}(\Omega)$ to $m^I\cdot\grad p[m^I]$.
Consequently, the limit passage $\eps\to 0$ in the identity \eqref{reg_Poisson_identity} gives
$\E_\eps(m^I)\to \E(m^I)$.
\endproof

To conclude the proof of Theorem~\ref{thm:global_ex}, we fix a $T>0$ and construct a a global solution $(m,p[m])$
on $(0,\infty)$ by concatenation of weak solutions on time intervals of length $T$ as constructed in Lemma~\ref{lem:eps_to_zero}.
This is possible due to the energy dissipation inequality \eqref{energy_ineq} and yields the global solution
announced in Theorem~\ref{thm:global_ex}.

\begin{remark}
Since the solution $(m,p)$ constructed in Theorem~\ref{thm:global_ex} satisfies
$m\cdot\grad p \in L^\infty(0,\infty; L^2(\Omega))$ and $\grad p \in L^\infty(0,\infty;L^2(\Omega))$,
implying $(m\cdot\grad p)\grad p \in L^\infty(0,\infty;L^1(\Omega))$ and $\partial_t m \in L^2((0,\infty)\times\Omega)$,
$|m|^{2\gamma-1} \in L^\infty(0,\infty;L^{2\gamma/(2\gamma-1))})$,
we conclude $\laplace m\in L^2(0,\infty;L^1(\Omega))$.
Theorems 1.7 and 3.3. in \cite{Guidetti} imply Besov regularity for $m$, namely $m\in L^2(0,\infty;B_\infty^{2,1}(\Omega))$.
Also, weak solutions satisfy the equation \eqref{eq1} pointwise almost everywhere
(while no regularity on second derivatives of $p$ is guaranteed).
Finally, we note that - by the same line of argument - weak stationary solutions posses the Besov regularity $m\in B_\infty^{2,1}(\Omega)$.
\end{remark}

\section{Existence and uniqueness of mild solutions}
\label{sec:mild_sols}

We fix $T>0$ and define the Banach spaces $\bbX:=(L^\infty(\Omega)\cap \mathrm{VMO}(\Omega))^d$ and ${\calX_T}:=L^\infty(0,T;\bbX)$,
where $\mathrm{VMO}(\Omega)$ denotes the space of functions with vanishing mean oscillation, see, e.g., \cite{Sarason}.
Note that the spaces $\bbX$ and, consequently, $\calX_T$, equipped with pointwise multiplication, are Banach algebras.
We denote $L:=D^2\laplace$, where $\laplace$ stands for the Dirichlet Laplacian on $\Omega$.
Moreover, we define the mapping $\calT$ on $\R\times\calX_T$ by
\( \label{T}
    \calT(\lambda, m) = e^{Lt}m^I + \int_0^t e^{L(t-s)} \bigl( \lambda F[m](s) - G[m](s) \bigr) \d s
\)
with $F[m] = (m\cdot\grad p[m])\grad p[m]$ where
$p[m]$ is the $H_0^1(\Omega)$-solution of the Poisson equation \eqref{eq1} with $m$ given,
and $G[m] = |m|^{2\gamma-1}m$.

Obviously, $(m,p)$ is a mild solution of the system \eqref{eq1}--\eqref{IC_0} with $\lambda=c^2$
subject to the initial datum $m^I$ if $m$ is a fixed point of $\calT$, i.e., $\calT(c^2,m)=m$.

The main result of this section is 
\begin{theorem} [Continuation from $c=0$]
Let $m^I\in \bbX$ and $S\in L^\infty(\Omega)$.
Then there exists an unbounded continuum of unique solutions $(\lambda, m)$ of $\calT(\lambda,m)=m$ in $[0,\infty)\times{\calX_T}$
emanating from $(0,m_0= e^{Lt}m^I)$.  
\label{thm:exMawhin}
\end{theorem}

To prove this theorem, we need a few auxiliary results.

\begin{lemma}\label{lem:ellReg}
Let $m\in\bbX$, $F\in L^q(\Omega)$ for some $1\leq q < \infty$ and $S\in L^r(\Omega)$ with $r=\max\{1,dq/(d+q)\}$.
Then the PDE
\[
    -\grad\cdot\bigl((I+m\otimes m)\grad p + F\bigr) &=& S\qquad\mbox{in }\Omega,\\
      p &=& 0\qquad\mbox{on }\partial\Omega
\]
has a unique weak solution $p\in W^{1,q}(\Omega)$
and there exists a constant $C(\Norm{m}_\bbX)>0$, independent of $F$ and $S$, such that
\(
    \Norm{\grad p}_{L^q(\Omega)} \leq C(\Norm{m}_{L^q(\Omega)}) \left( \Norm{F}_{L^q(\Omega)} + \Norm{S}_{L^r(\Omega)} \right).  \label{ellReg}
\)
\end{lemma}

\begproof
See Theorem 2.1 in \cite{Marino}.
\endproof

\begin{lemma}
Let $m^I\in \bbX$.
Then the Leray-Schauder degree of the mapping $I-\calT(0,\cdot)$, with $\calT$ given by \eqref{T}, satisfies
\[
   \mbox{deg}_{LS} [I-\calT(0,\cdot), {\calX_T}, 0] = \pm 1.
\]
\end{lemma}

\begproof
The maximum principle (applied componentwise) gives $m\in L^\infty((0,\infty)\times\Omega)^d$ for the unique
fixed point $m$ of $\calT(0,\cdot)$; see Lemma~\ref{lem:m-aux} for uniqueness.
Then, by standard parabolic theory we have
$m\in L^q(0,T; W^{2,q}(\Omega))^d \cap W^{1,q}(0,T; L^q(\Omega))^d$.
Since for $q$ large enough this space is embedded into ${\calX_T}$ (see Corollary 8 of \cite{Simon87}),
$m$ is the unique fixed point in $\calX_T$ of the mapping $\calT(0,\cdot)$.
\endproof

Our main existence result is based on Theorem 6.4 of \cite{Mawhin}, which we rephrase here
for the sake of reader's comfort:

\begin{theorem}[Theorem 6.4 of \cite{Mawhin}]\label{thm:Mawhin}
Let $\calT:[0,\infty)\times{\calX_T} \to \calX_T$ be a completely continuous mapping.
Assume that the Leray-Schauder degree $\mbox{deg}_{LS} [I-\calT(0,\cdot), {\calX_T} , 0]$
is well defined and non-zero.
Then there exists an unbounded continuum $\mathcal C$ of fixed points $\{(\lambda,m)\in[0,\infty)\times{\calX_T};\; \calT(\lambda,m)=m\}$
with $\mathcal{C} \cap (\{0\}\times{\calX_T}) \neq\emptyset$.
\end{theorem}

In fact, Theorem 6.4 of \cite{Mawhin} provides the existence of a pair of continua of $(\lambda,m)$,
one in $(0,\infty)\times{\calX_T}$ and the other one $(-\infty,0)\times{\calX_T}$.
In light of Remark~\ref{rem:non-pattern-formation}, only the continuum of solutions for $c^2 := \lambda>0$
is of interest in the context of network formation.

With the material developed so far we can conclude the proof of the main theorem of this section. 
\\

\noindent{\bf Proof of Theorem~\ref{thm:exMawhin}: }
We apply Theorem \ref{thm:Mawhin} to the mapping $\calT$ defined by \eqref{T}.
For this, we have to prove that $\calT$ is compact and continuous in $(0,\infty)\times{\calX_T}$.

To prove continuity of $\calT$, let us take a sequence $(\lambda^k, m^k)$ converging
as $k\to\infty$ to $(\lambda,m)$ in the norm topology of $(0,\infty)\times{\calX_T}$.
Lemma \ref{lem:ellReg} implies $\grad p[m]\in L^\infty(0,T; L^q(\Omega))$ for every $1\leq q < \infty$.
We denote $p^k := p[m^k]$ and calculate
\[
    -\grad\cdot\bigl( (I+m^k\otimes m^k)\grad(p^k-p[m])
       + (m^k\otimes m^k - m\otimes m)\grad p[m] \bigr) &=& 0 \qquad\mbox{in }\Omega,\\
       p^k - p[m] &=& 0 \qquad\mbox{on }\partial\Omega.
\]
Another application of Lemma \ref{lem:ellReg} gives strong convergence of $\grad p^k$ to $\grad p[m]$
in $L^\infty(0,T;L^q(\Omega))$, for every $1\leq q < \infty$.
Therefore, the sequence $\lambda^k F[m^k] = \lambda^k (m^k\cdot\grad p[m^k])\grad p[m^k]$
converges strongly in $L^\infty(0,T;L^q(\Omega))$ to $\lambda F[m] = \lambda (m\cdot\grad p[m])\grad p[m]$,
again for every $1\leq q < \infty$.
Convergence of the term $G[m^k]$ to $G[m]$ in $C([0,T]\times\Omega)$ is immediate.
Standard regularity results for the heat semigroup $e^{Lt}$ show then that $\calT(\lambda^k,m^k) \to \calT(\lambda,m)$
in $\calX_T$ as $k\to\infty$.

For compactness of $\calT$, we pick a weakly converging sequence $(\lambda^k, m^k)$ in $\R\times{\calX_T}$.
Due to Lemma~\ref{lem:ellReg}, the sequence $\grad p[m^k]$ is bounded in $L^q((0,T)\times\Omega)$
for any $q<\infty$, where $p[m^k]$ are the solutions of \eqref{eq1} with $m^k$.
Consequently, the term $\lambda^k F[m^k] = \lambda^k (m^k\cdot\grad p[m^k])\grad p[m^k]$ is bounded in 
$L^q((0,T)\times\Omega)$ for any $q<\infty$, which implies boundedness of $\bar m^k := \calT(\lambda^k,m^k)$
in $L^q(0,T;W^{2,q}(\Omega)$ and $\partial_t \bar m^k$ in $L^q((0,T)\times\Omega)$, again for any $q<\infty$.
Then, Corollary 8 in \cite{Simon87} implies that the sequence $\bar m^k$
is relatively compact in the space $C(0,T; W^{1,q}(\Omega))^d$.
Since $q$ can be chosen arbitrarily large, we have the continuous embedding
of $W^{1,q}(\Omega)^d$ into $\bbX$.
So we finally conclude that $\bar m^k$ is relatively compact in ${\calX_T}$.

The local uniqueness of solutions is obtained with a contraction argument.
Let us consider two mild solutions $(m^1,p^1)$, $(m^2,p^2)$ of the system \eqref{eq1}--\eqref{IC_0} and fix $q > 1$.
Taking the difference of the Poisson equations for $p_1$ and $p_2$ gives
\[
    \grad\cdot\Bigl( (I+m^1\otimes m^1)(\grad p^1-\grad p^2) \Bigr) =
       \grad\cdot \Bigl([(m^1-m^2)\otimes m^1]\grad p^2 + m^2\otimes(m^1-m^2)\grad p^2 \Bigr).
\]
We then have the estimate (Theorem 2.1 in \cite{Marino})
\[
    \Norm{\grad p^1-\grad p^2}_{L^\infty(0,T; L^q(\Omega))} &\leq&
      C(\Norm{m^1}_{\calX_T}) \left( \Norm{m^1}_{\calX_T} \Norm{m^1-m^2}_{\calX_T} \Norm{\grad p^2}_{L^\infty(0,T; L^q(\Omega))} \right. \\
      &&  \left. + \Norm{m^2}_{\calX_T} \Norm{m^1-m^2}_{\calX_T} \Norm{\grad p^2}_{L^\infty(0,T;L^q(\Omega))} \right) \\
    &=& C(\Norm{m^1}_{\calX_T}) \left( \Norm{m^1}_{\calX_T} + \Norm{m^2}_{\calX_T} \right) \Norm{\grad p^2}_{L^\infty(0,T;L^q(\Omega))} \Norm{m^1-m^2}_{\calX_T}.
\]
We define the mapping $\calH: {\calX_T}\to{\calX_T}$ by
\[
    \calH(m) = e^{Lt}m^I + \int_0^t e^{L(t-s)}\left(c^2 \grad p\otimes\grad p m - |m|^{2(\gamma-1)}m \right)(s) \d s,
\]
and set $\tilde m^i := \calH(m^i)$ for $i=1,2$. The difference $\tilde m^1-\tilde m^2$ reads
\[
   \tilde m^1-\tilde m^2 = \int_0^t e^{L(t-s)}  \Bigl[ \left( |m^2(s)|^{2(\gamma-1)}m^2(s) - |m^1(s)|^{2(\gamma-1)}m^1(s) \right) \\
       + c^2 \left( \grad p^1(s)\otimes\grad p^1(s)m^1(s) - \grad p^2(s)\otimes\grad p^2(s) m^2(s) \right)\Bigr].
\]
We write the last term as
\[
    c^2\left( (\grad p^1\otimes\grad p^1)m^1 - (\grad p^2\otimes\grad p^2)m^2 \right) &=&
       c^2 (\grad p^1\otimes\grad p^1)(m^1-m^2) \\ && + c^2\Bigl(\grad p^1\otimes\grad (p^1-p^2) + \grad (p^1-p^2)\otimes\grad p^2\Bigr)m^2
\]
and estimate
\[
    \Norm{ c^2\left( (\grad p^1\otimes\grad p^1)m^1 - (\grad p^2\otimes\grad p^2)m^2 \right) }_{L^\infty(0,T;L^q(\Omega))}
      \leq c^2 \Norm{\grad p^1}^2_{L^\infty(0,T;L^{2q}(\Omega))} \Norm{m^1-m^2}_{\calX_T} \\
      + c^2 \left(\Norm{\grad p^1}^q_{L^\infty(0,T;L^q(\Omega))} + \Norm{\grad p^2}^q_{L^\infty(0,T;L^q(\Omega))} \right)
           \Norm{\grad (p^1- p^2)}^q_{L^\infty(0,T;L^q(\Omega))} \Norm{m^2}_{\calX_T} \\
        \leq C(\Norm{m^1}_{\calX_T},\Norm{m^2}_{\calX_T}) \Norm{m^1-m^2}_{\calX_T}.
\]
Moreover, we have
\[
   \Norm{ |m^1|^{2(\gamma-1)}m^1 - |m^2|^{2(\gamma-1)}m^2 }_{C([0,T];L^q(\Omega))} \leq
      C(\Norm{m^1}^{2(\gamma-1)}_{\calX_T},\Norm{m^2}^{2(\gamma-1)}_{\calX_T}) \Norm{m^1-m^2}_{\calX_T}.
\]
The semigroup property of $e^{Lt}$ in $L^q(\Omega)^d$ provides the estimates
\[
   \Norm{\tilde m^1-\tilde m^2}_{C([0,T]; L^q(\Omega))} &\leq& C_1(\Norm{m^1}_{\calX_T},\Norm{m^2}_{\calX_T})\,T\, \Norm{m^1-m^2}_{\calX_T},\\
   \Norm{\tilde m^1-\tilde m^2}_{L^q(0,T; W^{2,q}(\Omega))} &\leq& C_2(\Norm{m^1}_{\calX_T},\Norm{m^2}_{\calX_T}) \Norm{m^1-m^2}_{\calX_T}.
\]
Moreover, we have
\[
    \partial_t (\tilde m^1-\tilde m^2) &=& D^2\laplace(\tilde m^1-\tilde m^2) + \left( - |m^1|^{2(\gamma-1)}m^1 + |m^2|^{2(\gamma-1)}m^2 \right) \\
        && + c^2\left( (\grad p^1\otimes\grad p^1)m^1 - (\grad p^2\otimes\grad p^2)m^2 \right),
\]
which yields
\[
    \Norm{\partial_t (\tilde m^1-\tilde m^2)}_{L^q((0,T)\times\Omega)} \leq C_3(\Norm{m^1}_{\calX_T},\Norm{m^2}_{\calX_T}) \Norm{m^1-m^2}_{\calX_T}.
\]
Combining the above estimates gives
\[
   \Norm{\tilde m^1-\tilde m^2}_{W^{1,q}((0,T)\times\Omega)} \leq C_4(\Norm{m^1}_{\calX_T},\Norm{m^2}_{\calX_T}) \Norm{m^1-m^2}_{\calX_T}. 
\]
Choosing some $1\leq q_1 < q$, the H\"older inequality gives
\[
   \Norm{\tilde m^1-\tilde m^2}_{L^{q_1}(0,T; W^{2,q_1}(\Omega))} &\leq& \tilde C_2(\Norm{m^1}_{\calX_T},\Norm{m^2}_{\calX_T}) \Norm{m^1-m^2}_{\calX_T} T^{1/q_1-1/q},\\
   \Norm{\partial_t (\tilde m^1-\tilde m^2)}_{L^{q_1}((0,T)\times\Omega)} &\leq& \tilde C_3(\Norm{m^1}_{\calX_T},\Norm{m^2}_{\calX_T}) \Norm{m^1-m^2}_{\calX_T} T^{1/q_1-1/q},
\]
which implies
\[
   \Norm{\tilde m^1-\tilde m^2}_{W^{1,q_1}((0,T)\times\Omega)} \leq \tilde C_4(\Norm{m^1}_{\calX_T},\Norm{m^2}_{\calX_T}) \Norm{m^1-m^2}_{\calX_T} T^{1/q_1-1/q}.
\]
Then, for $q$ and $q_1$ sufficiently large, we have by the Sobolev embedding of $W^{1,q_1}((0,T)\times\Omega)$ into $C^{0,\alpha}([0,T]\times\bar\Omega)$,
\(   \label{contraction_est}
   \Norm{\calH(m^1)-\calH(m^2)}_{\calX_T} =
   \Norm{\tilde m^1-\tilde m^2}_{\calX_T} \leq \tilde C(\Norm{m^1}_{\calX_T},\Norm{m^2}_{\calX_T}) \Norm{m^1-m^2}_{\calX_T} T^{1/q_1-1/q}.
\)
Moreover, we write
\[
   \Norm{\calH(m)}_{\calX_T} \leq \Norm{\calH(0)}_{\calX_T} + \Norm{\calH(m)-\calH(0)}_{\calX_T} \leq
      \Norm{e^{L t}m^I}_{\calX_T} + \tilde C(\Norm{m}_{\calX_T},0) \Norm{m}_{\calX_T} T^{1/q_1 - 1/q}.
\]
Now we choose
\[
   R > \sup_{0\leq t\leq 1} \Norm{e^{L t}m^I}_{\calX_T}
\]
and $0<T<1$ so small that
\[
   \tilde C(\alpha,0)\alpha T^{1/q_1 - 1/q} < R\qquad\mbox{for all } 0\leq\alpha\leq R.
\]
Then $\calH$ maps the set $\{m\in{\calX_T}; \Norm{m}_{\calX_T} \leq R\}$ into itself,
and due to \eqref{contraction_est}, it is a contraction on ${\calX_T}$.
This implies the local uniqueness of the above constructed fixed points of $\calT(\lambda,\cdot)$.
\endproof

\begin{remark}
The assertion of Theorem~\ref{thm:exMawhin} implies the following:
If for some $\lambda>0$ there is no fixed point of $\calT$ in ${\calX_T}$,
then there exists a \emph{bounded} sequence of $c_k^2=\lambda_k>0$
and a sequence of corresponding fixed points $m^k\in{\calX_T}$ of $T(\lambda_k,\cdot)$,
such that $\Norm{m^k}_{\calX_T} \to \infty$ as $k\to\infty$.
Moreover, the contraction argument implies that the fixed points $m$ of $\calT(c^2,m)$
are either global in time classical solutions of \eqref{eq1}--\eqref{IC_0},
or there exists a $t>0$ and a sequence $t_k\to t$ as $k\to\infty$
such that $\Norm{m}_{C([0,t_k];\bbX)}\to\infty$ as $k\to\infty$.
\end{remark}

\begin{remark}\label{rem:Meryers}
In the two-dimensional setting $d=2$ it is possible to apply Theorem \ref{thm:exMawhin}
in the space $L^\infty((0,T)\times\Omega)^d$ instead of $\calX_T$.
The estimate of Lemma \ref{lem:ellReg} is replaced by the Meyers estimate, Theorem 1 in \cite{Meyers},
which states that \eqref{ellReg} holds for \emph{some} $q>2$ if $m$ in bounded in $L^\infty((0,T)\times\Omega)^d$.
Then, due to the compact embedding of the space $W^{2,q}(\Omega)$ into $C^{0,\alpha}(\bar\Omega)$,
continuity and compactness of the mapping $\calT$ \eqref{calT} in the topology of $L^\infty((0,T)\times\Omega)^d$
is proven along the lines of the proof of Theorem \ref{thm:exMawhin}.
However, the contraction argument fails since the Sobolev embedding
of $W^{1,q_1}((0,T)\times\Omega)$ into $C^{0,\alpha}([0,T]\times\bar\Omega)$ only holds
for $q_1 > 3$.
\end{remark}

\begin{remark}\label{rem:globalMildEx_1D}
In the one-dimensional setting $d=1$ the branch of solutions constructed in Theorem \ref{thm:exMawhin}
is in fact global in $\lambda=c^2$ for every $T>0$. This follows from the $L^\infty$ bound on $\partial_x p$
provided by Lemma \ref{lem:auxPoisson1D} below.
Then, the maximum principle yields an a priori bound on $m$ in $\calX_T$ for every $T>0$.
In other words, a unique global in time mild solution exists for every value $\lambda=c^2$
and every $m^I\in L^\infty(0,1)$.
\end{remark}

\begin{lemma}\label{lem:auxPoisson1D}
Let $f\in L^1(0,1)$ and $b$ measurable on $(0,1)$ such that $b(x) \geq b_0 >0$ for all $x\in[0,1]$.
Let $p\in H^1_0(0,1)$ be the unique weak solution of
\(  \label{auxPoisson1D}
    -\partial_x \bigl(b(x)\partial_x p(x)\bigr) = f(x)
\)
on $[0,1]$ subject to the boundary conditions $p(0)=p(1)=0$.
\[
    |\partial_x p(x)| \leq \frac{2\Norm{f}_{L^1(0,1)}}{b(x)}\qquad\mbox{ for all } x\in(0,1) . 
\]
\end{lemma}

\begproof
We assume $b$ smooth enough and integrate \eqref{auxPoisson1D} on $(0,x)$,
\[
    b(x)\partial_x p(x) = -F(x) + B,
\]
where $F(x) = \int_0^x f(s)\d s$ and $B$ is an integration constant.
Dividing by $b(x)$ and integrating once again leads to
\[
    p(x) = - \int_0^x \frac{F(s)}{b(s)} \d s + B \int_0^x \frac{\d s}{b(s)}.
\]
The right boundary condition $p(1)=0$ gives the value for $B$,
\[
   B = \left(\int_0^1 \frac{F(s)}{b(s)} \d s\right) \left(\int_0^1 \frac{\d s}{b(s)}\right)^{-1},
\]
which immediately shows $|B| \leq \Norm{F}_{L^\infty(0,1)}$.
Using this in the above formula for $\partial_x p$ yields
\[
    |\partial_x p(x)| \leq \frac{|F(x)|}{b(x)} + \frac{|B|}{b(x)} \leq \frac{2\Norm{F}_{L^\infty(0,1)}}{b(x)}
       \leq \frac{2\Norm{f}_{L^1(0,1)}}{b(x)}
\]
and a density argument finishes the proof.
\endproof

\section{Long term convergence}
\label{sec:longterm}

Energy dissipation is not only useful for proving existence of solutions.
It also provides a powerful tool to prove a long term convergence result of the transient solutions towards steady states given by
\(
    -\grad\cdot(\grad p + (m\otimes m)\grad p) &=& S,  \label{eq1_stat}\\
    - D^2\laplace m - c^2(m\cdot\grad p)\grad p + |m|^{2(\gamma-1)}m &=& 0, \label{eq2_stat}
\)
subject to the homogeneous Dirichlet boundary conditions
\(   \label{BC_stat}
   m(t,x) = 0,\qquad p(t,x)=0 \qquad \mbox{for all } x\in\partial\Omega,\; t\geq 0.
\)
Our result about the long time behavior of the transient solution, Theorem \ref{thm:LongTimeBeh} below,
is based on the following slight modification of Lemma \ref{lem:weak-strong}.

\begin{lemma}\label{lem:weak-strong_L4}
Let $m^k$ be a sequence of functions converging to $m$ in the norm topology of $L^4(\Omega)$, and denote by $p^k\in H_0^1(\Omega)$
the corresponding weak solutions of \eqref{eq1}.
Then $p^k$ converges strongly in $H^1(\Omega)$ to $p$ as $k\to\infty$, the unique $H_0^1(\Omega)$-solution of \eqref{eq1}.
\end{lemma}

\begproof
We only need to slightly modify the proof of Lemma~\ref{lem:weak-strong}.
Due to the a priori bounds
\[
     \int_\Omega |\grad p^k|^2 \d x \leq C(S,\Omega),\qquad \int_\Omega |m^k\cdot\grad p^k|^2 \d x \leq C(S,\Omega)
\]
there exists a subsequence of $p^k$ converging weakly in $H^1(\Omega)$ to some $p\in H^1(\Omega)$.
Then, we can pass to the limit in the term $(m^k\otimes m^k)\grad p^k$ due to the assumed strong convergence
of $m^k$ in $L^4(\Omega)$.
We then continue along the lines of the proof of Lemma~\ref{lem:weak-strong}.
\endproof

\begin{theorem} [Long term convergence]
Let $(m,p)$ be a weak solution of the system \eqref{eq1}--\eqref{IC_0} constructed in Theorem~\ref{thm:global_ex}.
Fix $T > 0$ and a sequence $t_k \to \infty$ as $k\to\infty$
and define the time-shifts $m^{(t_k)}(\tau,x):=m(\tau+t_k,x)$, $p^{(t_k)}(\tau,x):=p(\tau+t_k,x)$ for $\tau\in (0,T)$.
Then there exists a subsequence, again denoted by $t_k$, such that, as $k\to\infty$,
\[
   m^{(t_k)}&\to& m^\infty\quad \mbox{strongly in } L^q(0,T; L^4(\Omega)) \mbox{ for any }q<\infty,\\
   p^{(t_k)}&\to& p^\infty\quad \mbox{strongly in } L^2(0,T; H^1_0(\Omega)),
\]
where $(m^\infty,p^\infty)$ only depends on $x$
and is a weak solution of the stationary system \eqref{eq1_stat}--\eqref{BC_stat}.
\label{thm:LongTimeBeh}
\end{theorem}

\begproof
The energy dissipation inequality \eqref{energy_ineq} implies
\[
    \int_0^\infty \int_\Omega |\partial_t m(t,x)|^2 \d x\d t < +\infty.
\]
Therefore,
\[
    \int_{t_k}^{t_k+T} \int_\Omega |\partial_t m(t,x)|^2 \d x\d t =
    \int_0^T \int_\Omega |\partial_\tau m^{(t_k)}(\tau,x)|^2 \d x\d\tau \to 0 \qquad\mbox{as } k\to\infty.
\]
Moreover, $\int_\Omega |\grad m|^2 \d x$ is uniformly bounded in time, so that
\[
    \int_0^T \int_\Omega |\grad m^{(t_k)}(\tau,x)|^2 \d x\d\tau < C
\]
for a constant $C$ independent of $k$.
Consequently, the sequence $m^{(t_k)}$ is uniformly bounded in $H^1((0,T)\times\Omega)$ and so
there exists a subsequence, again denoted by $t_k$, and $m^\infty\in H^1((0,T)\times\Omega)$ such that
$m^{(t_k)} \rightharpoonup m^\infty$ weakly in $H^1((0,T)\times\Omega)$.
Due to the weak lower semicontinuity of the norm, we have
\[
    \int_0^T \int_\Omega |\partial_\tau m^\infty(\tau,x)|^2 \d x\d\tau \leq
    \liminf_{k\to\infty} \int_0^T \int_\Omega |\partial_\tau m^{(t_k)}(\tau,x)|^2 \d x\d\tau = 0,
\]
so that $m^\infty$ is independent of $\tau$.

The Aubin-Lions compactness theorem yields the strong convergence of
$m^{(t_k)}$ to $m^\infty$ in the norm topology of $L^q(0,T; L^4(\Omega))$ for any $q<\infty$ if $d\leq 3$.
Therefore, a straightforward modification of Lemma~\ref{lem:weak-strong_L4}
yields strong convergence of $\grad p^{(t_k)}$ to $\grad p^\infty$ in $L^2((0,T)\times\Omega)$,
where $p^\infty$ is the solution of \eqref{eq1} with $m^\infty$.
Clearly, $p^\infty$ does not depend on $\tau$ as well.

It remains to show how to pass to the limit in the nonlinear terms of equation \eqref{eq2}.
Due to the energy dissipation inequality \eqref{energy_ineq}, the term
$m^{(t_k)}\cdot\grad p^{(t_k)}$ is uniformly bounded in $L^2((0,T)\times\Omega)$
and so it has a converging subsequence. The limit can be identified as $m^\infty\cdot\grad p^\infty$
due to the strong convergence of $m^{(t_k)}$ to $m^\infty$ in $L^2(0,T; L^4(\Omega))$
and of $\grad p^{(t_k)}$ to $\grad p^\infty$ in $L^2((0,T)\times\Omega)$.
By the same strong convergence, the whole term $(m^{(t_k)}\cdot\grad p^{(t_k)})\grad p^{(t_k)}$
converges to $(m^\infty\cdot \grad p^\infty)\grad p^\infty$ weakly in $L^1((0,T)\times\Omega)$.

Finally, the energy dissipation gives a uniform bound on $m^{(t_k)}$ in $L^\infty(0,T;L^{2\gamma}(\Omega))$,
and we pass to the limit in the algebraic term $|m^{(t_k)}|^{2(\gamma-1)}m^{(t_k)}$
due to Lemma~\ref{lem:m-alg-conv}.
\endproof

\section{The zero stationary state}
\label{sec:zero}

As shown in the previous section,  the stationary problem \eqref{eq1_stat}--\eqref{BC_stat} carries the information concerning  pattern formation. We present now several properties of this problem departing from the zero steady state. Indeed,  for all values of the parameters $D$ and $c$, the zero steady state is defined by  $m_0\equiv 0$ and $p_0\in H_0^1(\Omega)$ solving $-\laplace p_0 = S$ on $\Omega$. The main question that we are going to address here is whether nontrivial stationary solutions exist for certain parameter ranges.

One can try to prove existence of a non-zero stationary solution by variational methods. It is immediate to show that solutions of \eqref{eq1_stat}--\eqref{BC_stat}
are critical points of the functional
\[
    \widetilde E[m,p] = \frac12 \int_\Omega \left( D|\grad m|^2 + \frac{|m|^{2\gamma}}{\gamma}
       - c^2 |m\cdot\grad p|^2 - c^2|\grad p|^2 + 2c^2 Sp \right) \d x
\]
defined for $m\in H_0^1(\Omega)^d \cap L^{2\gamma}(\Omega)^d$, $p\in H_0^1(\Omega)^d$
such that $m\cdot\grad p \in L^2(\Omega)$.
Obviously, $\widetilde E$ does not posses any of the classical properties
that provide the existence of nontrivial critical points
(indeed, $\widetilde E$ is not bounded below, not convex, and does not
render to an application of the Mountain Pass Theorem in a straightforward way).
Another approach to stationary solutions is to consider $p=p[m]$,
the unique solution of the Poisson equation
\[
    -\grad\cdot[\grad p + (m\otimes m)\grad p] &=& S \qquad \mbox{in }\Omega,\\
    p &=& 0 \qquad \mbox{on }\partial\Omega.
\]
Note that by Lemma~\ref{lem:weak-strong_L4}, the unique solution $p=p[m]$ exists for any given $m\in L^4(\Omega)$.
Then one may try to find critical points of the energy functional \eqref{energy},
\[
   \E(m) = \frac12 \int_\Omega \left( D|\grad m|^2 + \frac{|m|^{2\gamma}}{\gamma} + c^2 |m\cdot\grad p[m]|^2 + c^2 |\grad p[m]|^2 \right) \d x,
\]
defined for $m\in H_0^1(\Omega)^d \cap L^{2\gamma}(\Omega)^d$.
Note that $\E(m) = \widetilde E[m,p[m]]$ and $D_p \widetilde E[m,p[m]]=0$,
so that critical points of $\widetilde E$ correspond to critical points of $\E$, and vice versa.

We shall later on resort to analyze bifurcations off zero steady state that we will call the \emph{branch of trivial stationary solutions} $(m_0\equiv 0,p_0)$. 
In order to see in which range of parameters this is possible, we first show a negative result.

\subsection{Stability of the zero steady state for $D$ large}

We show that for $D\to\infty$, weak solutions of the stationary system \eqref{eq1_stat}--\eqref{BC_stat} converge to zero. As a consequence, we cannot expect pattern formation when $D$ is too large.

\begin{proposition}
Let $(m^D, p^D)$ be a sequence of weak solutions of the stationary system \eqref{eq1_stat}--\eqref{BC_stat}
with the diffusion constant $D>0$. Then
\[
    m^D \to 0 \quad\mbox{in } H^1_0(\Omega),\qquad p^D\to p_0\quad\mbox{in } H^1_0(\Omega)\qquad \mbox{as } D\to\infty,
\]
where $p_0$ is the unique solution of
\(  \label{p0_problem}
   -\laplace p_0 = S \quad\mbox{in } \Omega, \qquad p_0=0 \quad\mbox{on } \partial\Omega.
\)
\end{proposition}

\begproof
We will skip the superscripts in $(m^D, p^D)$ for the sake of better legibility.
Multiplication of \eqref{eq1_stat} by $p$ and \eqref{eq2_stat} by $m$ and integration by parts yields
\[
    \int_\Omega |\grad p|^2 \d x + \int_\Omega |m\cdot\grad p|^2 \d x &=& \int_\Omega Sp,\\
    D^2 \int_\Omega |\grad m|^2 \d x - c^2 \int_\Omega |m\cdot\grad p|^2 \d x + \int_\Omega |m|^{2\gamma}\d x &=& 0.
\]
Multiplication of the first identity by $c^2$ and subtraction from the second gives
\[
   D^2 \int_\Omega |\grad m|^2 \d x + c^2 \int_\Omega |\grad p|^2 \d x + \int_\Omega |m|^{2\gamma}\d x = c^2 \int_\Omega Sp\d x.
\]
With an application of the Poincar\'e inequality, this implies
\[
   D^2 \int_\Omega |\grad m|^2 \d x + \frac{c^2}{2} \int_\Omega |\grad p|^2 \d x + \int_\Omega |m|^{2\gamma}\d x \leq C \int_\Omega S^2 \d x,
\]
and the strong convergence of $m$ to zero in $H^1_0(\Omega)$ follows.

Thanks to the Sobolev embedding, we also have strong the convergence $m\to 0$ in $L^4(\Omega)$ for $d\leq 3$,
and a slight modification of Lemma~\ref{lem:weak-strong} implies the strong convergence of $p$ to $p_0$ in $H^1_0(\Omega)$.
\endproof

A stronger result can be shown for the spatially one-dimensional case, namely,
that for large enough diffusivities $D$ the stationary problem only admits trivial solutions.

\begin{proposition}\label{prop:bigD_1D}
There exists a $D_0>0$ such that if $(m,p)$ is a solution of the 1D stationary problem
\(
   - \partial_x \bigl(\partial_x p + m^2 \partial_x p \bigr) &=& S,  \label{eq1D1_stat}\\
   - D^2\partial^2_{xx}m - c^2(\partial_x p)^2 m + |m|^{2(\gamma-1)}m &=& 0,  \label{eq1D2_stat}
\)
for $x\in(0,1)$, subject to homogeneous Dirichlet boundary conditions for $m$ and $p$ at $x\in\{0,1\}$ and
with the diffusion constant $D > D_0$, then $m\equiv 0$ almost everywhere on $(0,1)$
and $p$ is the weak solution of $-\partial_{xx}^2 p = S$ on $(0,1)$ with $p(0)=p(1)=0$.
\end{proposition}

\begproof
Assume that $\int_0^1 m^2 \d x > 0$. Multiplication of \eqref{eq1D2_stat} by $m$ and integration by parts yields
\[
   D^2 \int_0^1 (\partial_x m)^2 \d x &=& c^2 \int_0^1 (\partial_x p)^2 m^2 \d x - \int_0^1 |m|^{2\gamma} \d x \\
     &\leq& 4c^2 \Norm{S}^2_{L^1(0,1)} \int_0^1 m^2 \d x \\
     &\leq& 4C_1^2 c^2 \Norm{S}^2_{L^1(0,1)} \int_0^1 (\partial_x m)^2 \d x,
\]
where we used the uniform bound on $\partial_x p$ provided by Lemma~\ref{lem:auxPoisson1D} with $b(x)=1+m(x)^2$,
followed by the Poincar\'e inequality with constant $C_1>0$.
The above immediately implies that if $D>D_0$ for $D^2_0 := 4C_1^2 c^2 \Norm{S}^2_{L^\infty(0,1)}$,
then $m=0$ almost everywhere in $(0,1)$.
\endproof

\subsection{Bifurcations off the branch of trivial solutions}

We study the existence of nontrivial solutions of the stationary system \eqref{eq1_stat}--\eqref{BC_stat}
using a global bifurcation theorem by P. Rabinowitz \cite{Rabinowitz} and a local one for variational problems 
by the same author \cite{Rabinowitz-MPT}.
We will first consider the case $\gamma >1$ and assume that $S\not\equiv 0$ is smooth on $\Omega$.
We again work in the space $\bbX=(L^\infty(\Omega) \cap VMO(\Omega))^d$ here.

We decompose the solution of the Poisson equation \eqref{eq1_stat} as $p=p_0 + q[m]$,
where $p_0$ is the unique solution of \eqref{p0_problem} and $q=q[m]$ solves
\[
    -\grad\cdot[\grad q + (m\otimes m)\grad q] &=& \grad\cdot[(m\otimes m)\grad p_0]\quad\mbox{in }\Omega,\\
    q &=& 0\quad\mbox{on }\partial\Omega.
\]
The assumption $S\not\equiv 0$ implies $\grad p_0 \not\equiv 0$.
Let us fix $D>0$ and introduce the notations $\beta:=c^2/D^2$ and
\(
    Lm := (-\laplaceD)^{-1}(\grad p_0\otimes\grad p_0)m,  \label{L(m)}
\)
where $\laplaceD$ denotes the Dirichlet Laplacian on $\Omega$.
We also define the set
\[
    \calR(L) := \{\beta_0\in\R;\, \exists m_0\in\bbX,\,m_0\neq 0\mbox{ such that } m_0=\beta_0 Lm_0 \}.
\]
The stationary system \eqref{eq1_stat}--\eqref{BC_stat} is then equivalent to the fixed point problem
\(   \label{fixpoint_bif}
    m = \beta Lm + F(m,\beta).
\)
with $L$ given by \eqref{L(m)} and
\[
    F(m,\beta) &:=& \beta (-\laplaceD)^{-1}\Bigl( (\grad p_0\otimes\grad q + \grad q\otimes\grad p_0
       + \grad q\otimes\grad q)m - \frac{1}{D^2} |m|^{2(\gamma-1)}m \Bigr).
\]
Recall that for $\beta\leq 0$ the problem \eqref{fixpoint_bif} has only the trivial solution $m\equiv 0$,
as pointed out in Remark~\ref{rem:non-pattern-formation}. Therefore, we restrict the formulation of the following Theorem
to $\beta>0$.

\begin{proposition}
At every point $(m_0\equiv 0,\beta_0>0)\in\bbX\times\R$ for which $\beta_0\in\calR(L)$
there is a bifurcation off the branch of trivial solutions $(m\equiv 0,\beta)$
of a solution branch of the stationary system \eqref{fixpoint_bif}.
The branch either meets $\infty$ in $\bbX\times\R$ or meets a point
$(m_0\equiv 0,\beta_1)$ where $\beta_1\in\calR(L)$.
\end{proposition}

\begproof
Similarly as in the proof of Theorem~\ref{thm:exMawhin} it can be shown that the operators
$L$ and $F$ are continuous and compact as mappings from $\bbX$ to itself.
Due to the embedding of $C^{0,\alpha}(\Omega)$ into $VMO(\Omega)$, we have
\[
    \Norm{F(m,\beta)}_\bbX \leq C_\alpha \Norm{F(m,\beta)}_{C^{0,\alpha}(\bar\Omega)}\qquad\mbox{for all } 0 <\alpha\leq 1,
\]
and due to the Sobolev embedding
\[
   \Norm{F(m,\beta)}_\bbX \leq C_{\alpha,r} \Norm{F(m,\beta)}_{W^{2,r}(\Omega)}
\]
for $r$ sufficiently large.
Consequently, the estimate
\[
   \Norm{q[m]}_{W^{1,r}(\Omega)} \leq C_r(\Norm{m}_\bbX)\Norm{m}^2_{L^\infty(\Omega)}\qquad\mbox{for all } 1\leq r <\infty
\]
provided by \cite{Marino} gives
\[
   \Norm{F(m,\beta)}_\bbX \leq C_{\alpha,r}(\Norm{m}_\bbX,\beta)\left(\Norm{m}_{L^\infty(\Omega)}^3 + \Norm{m}^{2\gamma-1}_{L^\infty(\Omega)} \right),
\]
where $C_{\alpha,r}(\Norm{m}_\bbX,\beta)$ is bounded on bounded subsets of $\R^2$.

To apply the Theorem of Rabinowitz (Theorem 1.3 in \cite{Rabinowitz}), we need to study the eigenvalue problem
$n=\mu L(n)$, i.e.,
\[
    -\laplace n &=& \mu(\grad p_0\otimes\grad p_0)n \qquad \mbox{in }\Omega,\\
    n&=& 0 \qquad \mbox{on }\partial\Omega.
\]
Introducing the new variable $u:=(-\laplaceD)^{1/2}n$, the above problem is written as $u = \mu H u$, with
\(
    H u := (-\laplaceD)^{-1/2} (\grad p_0\otimes\grad p_0)(-\laplaceD)^{-1/2}u.   \label{calT}
\)
It is easy to prove that $H$ is a self-adjoint and compact operator $L^2(\Omega)^d \to L^2(\Omega)^d$,
and $\int_\Omega u\cdot Tu \d x \geq 0$ for all $u\in L^2(\Omega)^d$.
Consequently, the Spectral Theorem \cite{Brezis} implies that the spectrum of $H$ consists of a sequence of nonnegative real eigenvalues
$\sigma_1 > \sigma_2 > \dots \geq 0$ and possibly zero. Moreover,
\[
   \sigma_1 &=& \sup_{\Norm{u}_{L^2(\Omega)^d}=1} \int_\Omega |\grad p_0\cdot(-\laplaceD)^{-1/2} u|^2 \d x\\
   &=& \sup_{\Norm{v}_{H_0^1(\Omega)^d}=1} \int_\Omega |\grad p_0\cdot v|^2 \d x.
\]
Consequently, Theorem 1.3 of \cite{Rabinowitz} provides the branch of bifurcating nontrivial solutions
at eigenvalues of $H$ with odd multiplicity.
Bifurcation off eigenvalues with even multiplicity follows from the fact that stationary solutions $m$
are critical points of the energy functional $\E(m)$ \eqref{energy}, by applying the theory of \cite{Rabinowitz-MPT},
based on a local application of the mountain-pass theorem, combined with the global techniques of \cite{Rabinowitz}.
\endproof

\begin{remark}
The same result can be obtained for the case $\gamma=1$ by replacing $-\laplaceD$ by $-\laplaceD + \frac{1}{D^2}I$.
\end{remark}

Note that $\sigma_n\to 0$ as $n\to\infty$.
In the case $d=1$, zero is an eigenvalue if there exists a subinterval of $\Omega=(0,1)$ where $\partial_x p_0 = 0$.
For dimensions $d>1$, zero is always an eigenvalue of $T$,
since one can construct an eigenfunction $u\not\equiv 0$ such that $\grad p_0\cdot (-\laplaceD)^{-1/2}u = 0$.
Since $\mbox{dim}(T-\sigma_n I)<\infty$ for all $n\in\N$, and since the union of all
eigenspaces and $N(T)$ is the whole space $L^2(\Omega)^d$, we find that infinitely many eigenvalues $\sigma_n >0$ exist.

In the one-dimensional case $d=1$ it follows that the largest eigenvalue
$\sigma_1$ of $H$ is of odd multiplicity. Indeed, $\sigma_1$ is also the largest eigenvalue
of the mapping $(-\partial_{xx}^2)^{-1}(\partial_x p_0)^2$ with homogeneous Dirichlet boundary conditions.
Note that this mapping leaves the positive cone in $L^2(\Omega)$ invariant, so by the Hess-Kato extension
\cite{Hess-Kato} of the Krein-Rutman Theorem, the multiplicity of its largest eigenvalue is $1$
and the corresponding eigenfunction $n_1$ is nonnegative.
Therefore, we have the bifurcation of two continua of solutions from $(m_0\equiv 0, \beta)$
at $\beta=1/\sigma_1$, locally parametrized by
\[
    (m^\alpha,\beta^\alpha) = (\alpha n_1 + o(1), 1/\sigma_1 + o(1)) \qquad\mbox{for }\alpha\mbox{ close to }0.
\]
Due to the a priori estimates on $m$ in $\calX$ (see proof of Theorem~\ref{thm:exMawhin}),
both bifurcating branches either contain nontrivial solutions for all $\beta > 1/\sigma_1$
or they meet the trivial branch $(m_0\equiv 0,\beta)$ at another eigenvalue $\beta_1\in\calR(L)$.

\subsection{Linearized instability of the zero steady state}

We recall that  $\beta := c^2/D^2$ and work in any spatial dimension. We assume $\gamma>1$.

\begin{lemma}
Let $\sigma_1$ be the largest eigenvalue of $L$ defined in \eqref{L(m)}.
The trivial solution $(m_0,p_0)$ of the system \eqref{eq1}--\eqref{BC_0} with $\gamma>1$ is linearly asymptotically stable if
$\beta < 1/\sigma_1$, non-asymptotically stable if $\beta = 1/\sigma_1$ and exponentially unstable
if $\beta > 1/\sigma_1$.
\end{lemma}

This result explains why for $D$ small enough, patterns may occur as a consequence of  Turing instability. Typical is that the ratio of diffusions between the two equations should be correctly ordered and that a priori bounds on the steady state solutions exclude highly oscillatory instabilities. 
\\

\begproof
The linearization (G\^{a}teaux derivative) of \eqref{eq2} around $(m_0\equiv 0,p_0)$ in the direction $(n,q)$ reads
\[
   \partial_t n = D^2 M_\beta n := D^2(\laplace + \beta \grad p_0\otimes\grad p_0)n
\]
subject to $n(t=0,x)=n^I$ in $\Omega$ and homogeneous Dirichlet boundary conditions on $\partial\Omega$.
Let us define the quadratic form on $H^1_0(\Omega)^d \times H^1_0(\Omega)^d$,
\[
   Q_\beta(u,v) := (M_\beta u,v)_{L^2(\Omega)} = - \int_\Omega \grad u\cdot\grad v \d x + \beta \int_\Omega (\grad p_0\cdot u)(\grad p_0\cdot v) \d x.
\]
Then $Q_\beta$ is negative definite, i.e., $Q_\beta(n,n) <0$ for all $0\neq n\in H_0^1(\Omega)^d$, iff
\[
   \beta < \frac{\int_\Omega |\grad n|^2 \d x}{\int_\Omega |\grad p_0\cdot n|^2 \d x} \qquad\mbox{for all } 0\neq n\in H_0^1(\Omega)^d,
\]
i.e.,
\[
   \beta \leq \inf_{n\in H_0^1(\Omega)^d} \frac{\int_\Omega |\grad n|^2 \d x}{\int_\Omega |\grad p_0\cdot n|^2 \d x} = \frac{1}{\sigma_1},
\]
where the last equality is due to the formula for $\sigma_1$,
\[
     \sigma_1 = \sup_{n\in H_0^1(\Omega)^d} \frac{\int_\Omega |\grad p_0\cdot n|^2\d x}{\int_\Omega |\grad n|^2 \d x}.
\]
Consequently, $M_\beta$ has only negative eigenvalues and the trivial solution is linearly asymptotically stable if $\beta < 1/\sigma_1$.

On the other hand, there exists a nonnegative eigenvalue of $M_\beta$ iff $Q_\beta(n,n)\geq 0$ for some $0\neq n\in H_0^1(\Omega)^d$.
This happens if
\[
   \beta \geq \inf_{n\in H_0^1(\Omega)^d} \frac{\int_\Omega |\grad n|^2 \d x}{\int_\Omega |\grad p_0\cdot n|^2 \d x} = \frac{1}{\sigma_1}.
\]
Consequently, for $\beta=1/\sigma_1$, the largest eigenvalue of $M_\beta$ is zero and we have linearized non-asymptotic stability
(i.e., existence of a constant mode). Finally, for $\beta > 1/\sigma_1$ an exponentially growing mode exists.
\endproof

\section{Non-zero stationary states and pattern formation}
\label{sec:pattern formation}

As a consequence of previous results, we can expect non-zero steady states when $D$ is small. We can indeed build such steady states and analyze their stability. We proceed with the most general construction in one dimension

\subsection{One-dimensional case: nonlinear stability analysis with $D=0$}
\label{sec:stability1D}

We show now how stationary network-patterns are produced by the system in a special one dimensional setting on the interval $(0,1)$ with $D=0$ and we select those which are reachable by the dynamics.

The system \eqref{eq1}--\eqref{eq2} in one dimension with $D=0$ reads
\(
   - \partial_x \bigl(\partial_x p + m^2 \partial_x p \bigr) &=& S,  \label{eq1D1}\\
   \partial_t m - c^2(\partial_x p)^2 m + |m|^{2(\gamma-1)}m &=& 0,  \label{eq1D2}
\)
and, for the sake of simplicity, we consider it on the interval $(0,1)$ with mixed boundary conditions
\[
    \partial_x p(0) = 0, \qquad p(1) = 0\,.
\]
Integrating the first equation with respect to $x$ and taking into account the boundary conditions for $p$, we obtain
\[
    (1+m^2)\partial_x p = - \int_0^x S(y) \d y.
\]
Let us denote $B(x) := \int_0^x S(y)\d y \geq 0$, then we have
\(   \label{partx_p}
    \partial_x p = - \frac{B(x)}{1+m^2}.
\)
Inserting this into the equation for $m$ yields
\(   \label{ODE1}
    \partial_t m = \left( \frac{c^2 B(x)^2}{(1+m^2)^2} - |m|^{2(\gamma-1)} \right) m,
\)
which we interpret as a family of ODEs for $m=m(t)$ with the parameter $x$.

We now distinguish two cases:
\begin{itemize}
 \item $\gamma>1$: In this case the equation
 \(   \label{eq_m_s}
   \frac{c^2 B(x)^2}{(1+m^2)^2} - |m|^{2(\gamma-1)} = 0
 \)
 has, for a fixed $x\in(0,1)$, exactly two solutions $\pm m_s$ for some $m_s = m_s(x) >0$;
 assuming positivity of $S$, we have $B(x)^2 > 0$ on $(0,1)$.
 Thus, for every $x\in (0,1)$ the ODE \eqref{ODE1} has three stationary points, $m_0=0$, $m_s$ and $-m_s$.
 It can be easily checked that $m_0=0$ is unstable, while the other two are asymptotically stable.
 Therefore, solving \eqref{ODE1} subject to the initial datum $m^I=m^I(x)$ on $(0,1)$,
 we obtain, as $t\to\infty$, the asymptotic steady state $m_s(x)\sign(m^I(x))$.
 
 \item $\gamma=1$: In this case we have to solve the equation
 \[
   \frac{c^2 B(x)^2}{(1+m^2)^2} - 1 = 0,
 \]
 which has distinct nonzero real roots $\pm m_s(x)$ if and only if $cB(x) > 1$, with $m_s(x) = \sqrt{cB(x)-1}$.
 In this case, again, the ODE \eqref{ODE1} has three stationary points, unstable $m_0=0$ and stable $\pm m_s$.
 On the other hand, if $cB(x) \leq 1$, \eqref{ODE1} has the only stationary point $m=0$, which is stable.
 Thus, the solution of \eqref{ODE1} subject to the initial datum $m^I=m^I(x)$ on $(0,1)$
 converges to the asymptotic steady state $m_s(x)\sign(m^I(x))\chi_{\{cB(x) > 1\}}$.
\end{itemize}

\subsection{Linearized stability analysis of the system with $D=0$, $\gamma>1$}
\label{subsec:LinStabD=0}

We consider the system
\(
    -\div [(I + m\otimes m)\grad p] &=& S,   \label{D01}\\
    \part{m}{t} - c^2(m\cdot\grad p)\grad p + |m|^{2(\gamma-1)}m &=& 0, \label{D02}
\)
posed on a bounded domain $\Omega\subset\R^d$, subject to homogeneous Dirichlet boundary conditions
\[
    p = 0\,,\qquad m = 0 \qquad \mbox{for } x\in\partial\Omega, t\geq 0.
\]
For fixed $x$, there are  three kinds of stationary solutions of the problem \eqref{D02},
namely
\[
m_0(x) = \theta(x)\grad p_0(x) \qquad \text{with} \qquad    \theta(x) \in\left\{ 0, \pm c^\frac1{\gamma-1}|\grad p_0(x)|^\frac{2-\gamma}{\gamma-1} \right\}. 
\]
Once $\theta(x)$ has been chosen accordingly pointwise almost everywhere,
the stationary pressure $p_0$ satisfies
\[
   - \grad\cdot((I + m_0\otimes m_0) \grad p_0) = S.
\]
We construct the general stationary solution by fixing measurable disjoint sets
$\A_+\subseteq\Omega$, $\A_-\subseteq\Omega$, set $\A_0:=\Omega\setminus(\A_+\cup\A_-)$ and
\(  \label{m0}
    m_0(x) := \left(\chi_{\A_+}(x) - \chi_{\A_-}(x)\right) c^\frac1{\gamma-1}|\grad p_0(x)|^\frac{2-\gamma}{\gamma-1} \grad p_0(x),
\)
where $p_0$ solves
\(   \label{p0}
   -\grad\cdot\left[\left(1 + c^\frac{2}{\gamma-1}|\grad p_0(x)|^\frac{2}{\gamma-1}\chi_{\A_+\cup\A_-}(x) \right)\grad p_0(x) \right] = S. 
\)

\begin{theorem} [Existence of network-patterns]
For any $S\in L^2(\Omega)$, $\gamma>1$
and for any pair of measurable disjoint sets $\A_+$, $\A_-\subseteq\Omega$
there exists a unique weak solution $p_0\in H^1_0(\Omega) \cap W_0^{1,2\gamma/(\gamma-1)}(\A_+\cup\A_-)$ of \eqref{p0}.
\label{thm:exStat_gamma>1}\end{theorem}

\begproof
We set  $\A:=\A_+\cup\A_-$ and define the functional $\calF: H^1_0(\Omega) \to \R\cup \{+\infty\}$,
\[
   \calF[p] := \frac12 \int_\Omega |\grad p|^2 \d x + c^\frac{2}{\gamma-1} \frac{\gamma-1}{2\gamma} \int_\A |\grad p|^{\frac{2\gamma}{\gamma-1}} \d x
     - \int_\Omega p S \d x,
\]
and $\calF[p] := \infty$ if $\grad p\notin L^{\frac{2\gamma}{\gamma-1}}(\Omega)$.
Then $\calF$ is uniformly convex since $\frac{2\gamma}{\gamma-1} > 2$. Also, coercivity on $H^1_0(\Omega)$
is standard.
The classical theory (see, e.g., \cite{Evans}) provides then the existence of a unique minimizer $p_0\in H^1_0(\Omega)$ of $\calF$,
which is the unique solution to the corresponding Euler-Lagrange equation \eqref{p0}.
\endproof

The linearization (G\^{a}teaux derivative) of \eqref{D01}, \eqref{D02} around $(m_0,p_0)$ in direction $(n,q)$ is given by
\(
    0 &=& -\grad\cdot\bigl[ (I+m_0\otimes m_0)\grad q + (m_0\otimes n + n\otimes m_0)\grad p_0 \bigr],  \label{D0lin1}\\ 
    \partial_t n &=& c^2\bigl[ (n\cdot\grad p_0)\grad p_0 + (m_0\cdot\grad q)\grad p_0 + (m_0\cdot\grad p_0)\grad q\bigr] \nonumber\\
      && - |m_0|^{2(\gamma-1)}n - 2(\gamma-1)|m_0|^{2(\gamma-2)} (m_0\cdot n) m_0, \label{D0lin2}
\)
subject to the homogeneous Dirichlet boundary conditions
\(
   q = 0,\quad n=0\qquad \mbox{on }\partial\Omega  \label{D0lin_BC}
\)
and the initial condition for $n$,
\(
   n(t=0) = n^I \qquad\mbox{in }\Omega. \label{D0lin_IC}
\)

\begin{theorem}[Linear stability of the  network-pattern]
Let $\A_+$ and $\A_-$ be such that $\mathrm{meas}(\A_+\cup\A_-) = \mathrm{meas}(\Omega)$.
Moreover, assume
$\{x\in\Omega;\, \grad p_0(x)=0\} \subseteq \{x\in\Omega;\, n^I(x)=0\}$.
Let $(q,n)$ be the solution of \eqref{D0lin1}--\eqref{D0lin_IC}.
Then
\[
    \lim_{t\to\infty} \int_\Omega |\grad q(t,x)|^2 \d x = 0,\qquad
    \lim_{t\to\infty} \int_\Omega |n(t,x)|^2 \d x = 0.
\]
\label{thm:LinAsStab} 
\end{theorem}

\begproof
We first establish that $\grad q(t,x)$ lies in $L^2(\Omega\times(0,\infty)$.
Multiplication of \eqref{D0lin2} by $n$ and integration by parts yields
\[
   \frac12\tot{}{t} \int_\Omega n^2\d x &=& c^2 \int_\Omega (n\cdot\grad p_0)^2 + (m_0\cdot\grad q)(n\cdot\grad p_0) + (m_0\cdot\grad p_0)(n\cdot\grad q) \d x \\
      &-& \int_\Omega |m_0|^{2(\gamma-1)}|n|^2 \d x - 2(\gamma-1) \int_\Omega |m_0|^{2(\gamma-1)}(m_0\cdot n)^2\d x \,.
\]
We have the identities
\( \label{identities}
    |m_0|^{2(\gamma-1)} = c^2 |\grad p_0|^2 \chi_{\A_+\cup\A_-} \qquad\mbox{and}\qquad
    |m_0|^{2(\gamma-2)} (m_0\cdot n)^2 = c^2 (\grad p_0\cdot n)^2 \chi_{\A_+\cup\A_-}.
\)
Moreover, multiplication of \eqref{D0lin1} by $q$ and integration by parts gives
\(   \label{LinPoissonEnergy}
    \int_\Omega (m_0\cdot\grad q)(n\cdot\grad p_0) + (m_0\cdot\grad p_0)(n\cdot\grad q)\d x
    = -\int_\Omega |\grad q|^2 + |m_0\cdot\grad q|^2 \d x,
\)
so that we have
\[
   \frac12\tot{}{t} \int_\Omega n^2\d x &=& c^2 \int_\Omega |n\cdot\grad p_0|^2 \d x
     - c^2 \int_\Omega \left( |\grad q|^2 + |m_0\cdot\grad q|^2 \right) \d x \\
      &-& c^2 \int_{\A_+\cup\A_-} |\grad p_0|^2 |n|^2 \d x - 2(\gamma-1) c^2 \int_{\A_+\cup\A_-} |\grad p_0\cdot n|^2\d x \,.
\]
According to the assumption $\mathrm{meas}(\A_+\cup\A_-)=\mathrm{meas}(\Omega)$, and with the Cauchy-Schwarz inequality
$|n\cdot\grad p_0|^2\leq |n|^2 |\grad p_0|^2$, we have due to $\gamma>1$,
\[
   \frac12\tot{}{t} \int_\Omega n^2\d x + c^2 \int_\Omega |\grad q|^2 \d x \leq 0.
\]
This implies $\grad q \in L^2(\Omega\times(0,\infty)$,
\(    \label{integrability_grad_q}
    \int_0^\infty \int_\Omega |\grad q(t,x)|^2 \d x\d t < +\infty.
\)

Using the assumption $\mathrm{meas}(\A_0)=0$ and the identities \eqref{identities} in \eqref{D0lin2} yields
\[
   \part{n}{t} = c^2\bigl( (3-2\gamma)(\grad p_0\cdot n)\grad p_0 - |\grad p_0|^2 n\bigr)
               + c^2 \bigl( (m_0\cdot\grad q)\grad p_0 + (m_0\cdot\grad p_0)\grad q \bigr),
\]
so that we can write
\[
    \part{n}{t} = A(x)n + f(t,x)
\]
with
\[
    A(x) &=& c^2(3-2\gamma)\grad p_0\otimes\grad p_0 - c^2|\grad p_0|^2 I,\\
    f(t,x) &=& c^2 \bigl( (m_0\cdot\grad q)\grad p_0 + (m_0\cdot\grad p_0)\grad q \bigr).
\]
The linear ODE initial value problem has the solution
\[
    n(t,x) = n_\mathrm{hom}(t,x) + n_p(t,x) := e^{A(x)t}n^I(x) + \int_0^t e^{A(x)(t-s)}f(s,x)\d s.
\]
It can be easily calculated that $A(x)$ has the $(d-1)$-fold eigenvalue $-c^2|\grad p_0|^2$
and the 1-fold eigenvalue $-2c^2 (\gamma-1)|\grad p_0|^2$.
Therefore, denoting $\alpha:=\min(1,2(\gamma-1))>0$, we have the estimate
\[
    |n_\mathrm{hom}(t,x)| \leq \exp\left(-c^2\alpha t |\grad p_0(x)|^2\right) |n^I(x)|,
\]
so that, recalling the assumption $\{x\in\Omega;\; \grad p_0(x) = 0\} \subseteq \{x\in\Omega;\; n^I(x)=0\}$,
\[
    \lim_{t\to\infty}  \int_\Omega |n_\mathrm{hom}(t,x)|^2 \d x = 0.
\]

It remains to estimate the inhomogeneous part $n_p$ of $n$.
We start with
\[
   |n_p(t,x)| &\leq& \int_0^t \exp\left(-c^2\alpha(t-s)|\grad p_0(x)|^2\right)|f(s,x)|\d s,\\
   |f(t,x)| &\leq& 2\theta |\grad p_0(x)|^2 |\grad q(t,x)| = 2c^\frac{1}{\gamma-1} |\grad p_0(x)|^\frac{\gamma}{\gamma-1} |\grad q(t,x)|,
\]
which gives
\[
    \int_\Omega |n_p(t,x)|^2 \d x \leq 4c^\frac{2}{\gamma-1} \int_\Omega |\grad p_0(x)|^\frac{2\gamma}{\gamma-1}
       \left( \int_0^t \exp\left(-c^2\alpha(t-s)|\grad p_0(x)|^2\right) |\grad q(s,x)| \d s\right)^2 \d x.
\]
Denoting
\[
   I(t,x) := |\grad p_0(x)|^\frac{2\gamma}{\gamma-1}
       \left( \int_0^t \exp\left(-c^2\alpha(t-s)|\grad p_0(x)|^2\right) |\grad q(s,x)| \d s\right)^2,
\]
the Cauchy-Schwarz inequality gives
\[
   I(t,x) \leq |\grad p_0(x)|^\frac{2\gamma}{\gamma-1} \int_0^t \exp\left(-{c^2\alpha}(t-s)|\grad p_0(x)|^2\right) \d s\;
       \int_0^t \exp\left(-{c^2\alpha}(t-s)|\grad p_0(x)|^2\right) |\grad q(s,x)|^2 \d s,
\]
and denoting $K:= \int_0^\infty \exp\left(-\frac{c^2\alpha}{2}s\right) \d s < \infty$, we have
\(   \label{bound_I}
   I(t,x) \leq K |\grad p_0(x)|^\frac{2\gamma}{\gamma-1} \int_0^t \exp\left(-{c^2\alpha}(t-s)|\grad p_0(x)|^2\right) |\grad q(s,x)|^2 \d s.
\)
Now, \eqref{integrability_grad_q} states that $|\grad q(s,x)|^2\in L^1(0,\infty)$ for almost every $x\in\Omega$.
Clearly, for those $x\in\Omega$ where $|\grad p_0(x)|>0$, we have
\[
    \lim_{t\to\infty} \chi_{[0,t]}(s) \exp\left(-{c^2\alpha}(t-s)|\grad p_0(x)|^2\right) = 0\qquad\mbox{for all } s>0.
\]
We shall employ the following technical Lemma:

\begin{lemma}\label{lem:technical}
Let $h:(0,\infty)\to\R$ be in $L^1(0,\infty)$ and $\sigma>0$.
Then
\[
   \lim_{t\to\infty} \int_0^t \exp(-\sigma(t-s)) h(s) \d s = 0.
\]
\end{lemma}

\begproof
Clearly, $\lim_{t\to\infty} \exp(-\sigma(t-s)) h(s) = 0$ for every $s\in(0,\infty)$.
Moreover, we have
\[
    \int_0^t \exp(-\sigma(t-s)) h(s) \d s = \int_0^\infty \chi_{(0,t)}(s) \exp(-\sigma(t-s)) h(s) \d s,
\]
and the result follows by an application of the Lebesgue dominated convergence theorem with the integrable majorant $h(s)$.
\endproof

\noindent We apply the above Lemma with $h(s):=|\grad q(s,x)|^2$ for every $x\in\Omega$ where $|\grad p_0(x)|>0$ and conclude
\[
   && \lim_{t\to\infty} \int_0^t \exp\left(-{c^2\alpha}(t-s)|\grad p_0(x)|^2\right) |\grad q(s,x)|^2 \d s \\ &=&
   \lim_{t\to\infty} \int_0^\infty \chi_{[0,t]}(s) \exp\left(-{c^2\alpha}(t-s)|\grad p_0(x)|^2\right) |\grad q(s,x)|^2 \d s
   = 0.
\]
Therefore, $\lim_{t\to\infty} I(t,x) = 0$ for almost all $x\in\Omega$.
Moreover, \eqref{bound_I} implies
\[
   I(t,x) \leq K |\grad p_0(x)|^\frac{2}{\gamma-1} \int_0^\infty |\grad q(s,x)|^2 \d s,
\]
so that, with \eqref{integrability_grad_q} and Theorem~\ref{thm:LongTimeBeh}, $I(t,x)$
is bounded by an integrable majorant in $x$ for all $t\geq 0$
(since $\Omega\setminus(\A_+\cup\A_-)$ has Lebesgue measure zero).
An application of the Lebesgue theorem yields then
\[
    \lim_{t\to\infty}  \int_\Omega |n_p(t,x)|^2 \d x = 0.
\]
The second part of the claim,
\[
   \lim_{t\to\infty} \int_\Omega |\grad q(t,x)|^2 \d x = 0,
\]
follows directly from the limit passage in \eqref{LinPoissonEnergy}.
\endproof

\begin{remark}
On the set $\{x\in\Omega;\, \grad p_0(x)=0\}$ equation \eqref{D0lin2} reduces to
$\part{n}{t} = 0$. Therefore, the restriction on the perturbation $\{x\in\Omega;\, \grad p_0(x)=0\} \subseteq \{x\in\Omega;\, n^I(x)=0\}$
of Theorem~\ref{thm:LinAsStab} is necessary for linearized asymptotic stability.

The assumption $\mathrm{meas}(\A_0)=0$ is necessary as well, since for $x\in\A_0$ equation \eqref{D0lin2} reduces to
\[
    \part{n}{t} = c^2 (n\cdot\grad p_0)\grad p_0,
\]
so that there are non-decaying modes.
\end{remark}

\begin{remark}
Let $\gamma = 1$. Then for each stationary state $(m_0,p_0)$ of \eqref{D01},\eqref{D02} there exists a measurable real valued
function $\lambda = \lambda(x)$ such that
\[
    m_0(x) = \lambda(x) \chi_{\{c|p_0|=1\}}(x) \grad p_0(x)
\]
and $p_0$ solves the highly nonlinear Poisson equation
\[
    -\grad\cdot\left[ \left(1+\frac{\lambda(x)^2}{c^2}\chi_{\{c|p_0|=1\}}(x) \right) \grad p_0 \right] = S
\]
subject to the homogeneous Dirichlet boundary condition $p_0=0$ on $\partial\Omega$.
\end{remark}

\begin{remark}
The results obtained in this Section need to be related to the nonlinear stability analysis
of Section~\ref{sec:stability1D} in the 1D case $\Omega=(0,1)$.
In particular, in the case $\gamma>1$, we concluded that the solution
of the 1D system converges to $m_s(x)\sign(m^I(x))$ as $t\to\infty$, with $m^I$ the initial condition for $m$,
where $m_s(x)$ is the positive solution of \eqref{eq_m_s}.
Fixing $m^I$ with $m^I(x)\neq 0$ for almost all $x\in\Omega$, we may set $m_0(x) := m_s(x)\sign(m^I(x))$,
which is obviously a stationary solution of the nonlinear 1D system \eqref{eq1D1}--\eqref{eq1D2}.
We assume $S>0$ on $\Omega$, so that $B(x):=\int_0^x S(y)\d y$ is a positive and increasing function.
Taking into account \eqref{partx_p}, we have $\mathrm{meas}\{\partial_x p=0\}=0$
and $\mathrm{meas}\{m_0=0\}=0$, so that the assumptions of Theorem~\ref{thm:LinAsStab} are satisfied.
The Theorem then states that solutions $(n,q)$ of the linearized system \eqref{D0lin1}--\eqref{D0lin2},
that we interpret as perturbations of the stationary solution $(m_0,p_0)$,
converge to zero in the $L^2$-sense.
However, if we perturb $m_0$, say with $m_0 + \eps\bar m$ with small $\eps>0$,
such that $\mathrm{meas}\{\sign(m_0) \neq \sign(m_0 + \eps\bar m)\}>0$
and feed this perturbed state as an initial datum for the nonlinear system \eqref{eq1D1}--\eqref{eq1D2},
then clearly the corresponding solution will converge to a different state than $m_0$ as $t\to\infty$.
This apparent contradiction is explained as follows:
The linearization \eqref{D0lin1}--\eqref{D0lin2} can be identified as the $L^2(\Omega)^d$-G\^{a}teaux
derivative of the nonlinear problem \eqref{D01}--\eqref{D02} but \emph{not} as the $L^2(\Omega)^d$-Fr\'echet derivative.
Thus linearized asymptotic $L^2(\Omega)^d$-stability does not imply nonlinear asymptotic stability.
Also, nonlinear asymptotic stability of the ODEs \eqref{ODE1} (for each $x$ fixed)
does not imply asymptotic stability of the PDE system \eqref{eq1D1}--\eqref{eq1D2},
albeit it implies that the PDE solution from an $L^2$-close initial state converges to an $L^2$-close steady state.

The linearized stability analysis together with the nonlinear one-dimensional analysis explain how patterns are formed
for small values of the diffusion constant $D$ and large times $t$.
Stable large time patterns for $m$ have the form \eqref{m0} (with smoothing for $D>0$ but small),
where the sets $\A_+$ and $\A_-$ are determined by the initial datum $m(t=0)=m^I$ and, due to small diffusion,
the set $\{m(t)=0\}$ has zero Lebesgue measure for $t>0$.
Clearly stability has to be interpreted in the sense that close-by initial data $m^I$ generate patterns
close to the one generated by the unperturbed initial datum.
\end{remark}

\section{The limit $D\to 0$ in the one dimensional setting}
\label{sec:Dtozero1D}

We consider the one dimensional network formation system
\(
   - \partial_x \bigl((1+m^2)\partial_x p \bigr) &=& S,  \label{eq1D_par_1}\\
   \partial_t m - D^2\partial^2_{xx}m - c^2(\partial_x p)^2 m + |m|^{2(\gamma-1)}m &=& 0,  \label{eq1D_par_2}
\)
on $\Omega=(0,1)$, with $S\in L^\infty(0,1)$, subject to the homogeneous Dirichlet boundary conditions
\(   \label{bc1D_par}
   m(0,t)=m(1,t)=0, \qquad p(0,t)=p(1,t)=0,\qquad \mbox{for } t\geq 0.
\)
We prove that the solution of \eqref{eq1D_par_1}--\eqref{eq1D_par_2} converges to zero for large $D$.
We first establish the following uniform in $D$ a priori estimate:
\begin{lemma}[Uniform BV bound]
Let $(m,p)\in C^1([0,T];C^2[0,1])\times C([0,T];C^2[0,1])$ be a classical solution of the system
\eqref{eq1D_par_1}--\eqref{eq1D_par_2},
subject to the boundary conditions \eqref{bc1D_par}.
Then there exists a constant $C=C(T,\Norm{S}_{L^\infty(0,1)})$ independent of $D$
such that
\[
   \max_{t\in[0,T]} \int_0^1 |\partial_x m| \d x < C.
\]
\label{lem:BV1D}
\end{lemma}

\begproof
We take the derivative of \eqref{eq1D_par_2} with respect to $x$,
\(  \label{eq_der_x}
   \partial^2_{tx} m - D^2\partial^3_{xxx}m - c^2(\partial_x p)^2 \partial_x m
      - 2c^2 (\partial^2_{xx} p)(\partial_x p)m + (2\gamma-1)|m|^{2(\gamma-1)}\partial_x m = 0.
\)
We multiply \eqref{eq_der_x} by $\sign(\partial_x m)$ and integrate over $\Omega=(0,1)$,
\[
    \tot{}{t} \int_0^1 |\partial_x m| \d x &=&
      D^2 \int_0^1 (\partial^3_{xxx} m)\,\sign(\partial_x m) \d x
      + c^2 \int_0^1 |\partial_x m|(\partial_x p)^2 \d x  \\
      && + 2c^2 \int_0^1 m\,\sign(\partial_x m) (\partial^2_{xx} p)\partial_x p \d x
      - (2\gamma-1) \int_0^1 |m|^{2(\gamma-1)} |\partial_x m| \d x.
\]
The Kato inequality \cite{Kato} for the first term of the right-hand side yields
\[
   D^2 \int_0^1 (\partial^3_{xxx} m)\,\sign(\partial_x m) \d x \leq
    D^2 \int_0^1 \partial^2_{xx} |\partial_x m| \d x =
    D^2 \Bigl[ \partial_x |\partial_x m| \Bigr]_{x=0}^1 =
    D^2 \Bigl[ (\partial^2_{xx} m)\,\sign(m) \Bigr]_{x=0}^1.
\]
Inserting the homogeneous Dirichlet boundary conditions for $m$
into \eqref{eq1D_par_2} yields $\partial^2_{xx}m(0,t)=\partial^2_{xx}m(1,t)=0$,
so that the above boundary term vanishes.
Consequently, we have
\[
   D^2 \int_0^1 (\partial^3_{xxx} m)\,\sign(\partial_x m) \d x \leq 0.
\]
The product rule for the Poisson equation \eqref{eq1D_par_1} yields
\[
   \partial^2_{xx}p = -\frac{S}{1+m^2} - \frac{2m(\partial_x m)(\partial_x p)}{1+m^2},
\]
so that
\[
   2c^2 \int_0^1 m\,\sign(\partial_x m) (\partial^2_{xx} p)\partial_x p \d x
   &=& - 2c^2 \int_0^1 \frac{m\,\sign(\partial_x m)}{1+m^2} (\partial_x p)S \d x
    - 4c^2 \int_0^1 \frac{m^2 |\partial_x m| (\partial_x p)^2}{1+m^2} \d x \\
   &\leq&  2c^2 \Norm{S}_{L^2(0,1)}\left( \int_0^1 |m|^2 (\partial_x p)^2 \d x \right)^{1/2}
    - 4c^2 \int_0^1 \frac{m^2 |\partial_x m| (\partial_x p)^2}{1+m^2} \d x.
\]
Multiplying \eqref{eq1D_par_1} by $p$, integrating by parts
and using Cauchy-Schwarz and Poincar\'e inequalities gives
\[
    \int_0^1 (\partial_x p)^2 \d x + \int_0^1 m^2 (\partial_x p)^2 \d x &=& \int_0^1 pS \d x \\
      &\leq& \Norm{S}_{L^2(0,1)}\Norm{p}_{L^2(0,1)}   \\
            &\leq& C \Norm{S}^2_{L^2(0,1)} + \Norm{\partial_x p}^2_{L^2(0,1)},
\]
for suitable constant $C > 0$, so that
\[
    \int_0^1 m^2 (\partial_x p)^2 \d x \leq C \Norm{S}^2_{L^2(0,1)},
\]
and
\[
    2c^2 \Norm{S}_{L^2(0,1)}\left( \int_0^1 |m|^2 (\partial_x p)^2 \d x \right)^{1/2} \leq C\Norm{S}_{L^2(0,1)}^2.
\]
We are left with the terms
\[
   c^2 \int_0^1 |\partial_x m|(\partial_x p)^2 \d x - 4c^2 \int_0^1 \frac{m^2 |\partial_x m| (\partial_x p)^2}{1+m^2} \d x
       &=& c^2 \int_0^1 |\partial_x m| (\partial_x p)^2 \left( 1 - \frac{4m^2}{1+m^2}\right) \\
   &\leq& c^2 \int_0^1 |\partial_x m| (\partial_x p)^2 \d x.   
\]
Now, an application of Lemma~\ref{lem:auxPoisson1D} with $b(x):=1+m(x)^2$ gives the estimate
$|\partial_x p| \leq \frac{2\Norm{S}_{L^1(0,1)}}{1+m(x)^2} \leq C$, so that the above expression is estimated from above by
$C \int_0^1 |\partial_x m| \d x$.
Altogether, we have
\[
   \tot{}{t} \int_0^1 |\partial_x m| \d x \leq C \left( \int_0^1 |\partial_x m| \d x +  \Norm{S}^2_{L^2(0,1)} \right).
\]
An application of the Gronwall lemma gives
\[
    \max_{t\in[0,T]} \int_0^1 |\partial_x m| \d x \leq C(T,\Norm{S}^2_{L^2(0,1)}) < \infty
\]
for every $T>0$.
\endproof

\begin{theorem}[Limit of vanishing diffusion]
Let $(m^D,p^D)\in C^1([0,T];C^2[0,1])\times C([0,T];C^2[0,1])$ be classical solutions of the system
\eqref{eq1D_par_1}--\eqref{bc1D_par} with the diffusion constant $D>0$ on the time interval $[0,T]$.
Then there exists a subsequence such that $(m^D,p^D)\to (m,p)$ as $D\to 0$, where $(m,p)$ is a solution of
\(
   - \partial_x \bigl((1+m^2)\partial_x p \bigr) &=& S,  \label{eq1D_D0_1}\\
   \partial_t m - c^2(\partial_x p)^2 m + |m|^{2(\gamma-1)}m &=& 0,  \label{eq1D_D0_2}
\)
subject to the homogeneous Dirichlet boundary conditions
\(   \label{bc1D_D0}
   m(t,0)=m(t,1)=0, \qquad p(t,0)=p(t,1)=0,\qquad \mbox{for } t\geq 0.
\)
\label{thm:Dtozero1D}
\end{theorem}

\begproof
According to Lemma~\ref{lem:BV1D}, the family $m^D$ is uniformly bounded in $L^\infty(0,T; BV[0,1])$.
Moreover, the energy dissipation given by Lemma~\ref{lem:energy} provides a uniform bound
on $\partial_t m^D$ in $L^2((0,T)\times(0,1))$.
Corollary 4 in \cite{Simon87} implies then strong convergence of a subsequence of $m^D$ to $m$ in $L^q((0,T)\times (0,1))$ for any $q<\infty$,
and Lemma~\ref{lem:weak-strong_L4} gives strong convergence of $\partial_x p^D$ to $\partial_x p$ in $L^2((0,T)\times(0,1))$.
This allows us to pass to the limit $D\to 0$ in the Poisson equation \eqref{eq1D_D0_1}.

Moreover, we have a uniform bound on $\partial_x p^D$ in $L^\infty((0,T)\times(0,1))$ by Lemma~\ref{lem:auxPoisson1D}.
Consequently, due to the strong convergence of $\partial_x p^D$ to $\partial_x p$ in $L^2((0,T)\times(0,1))$,
$(\partial_x p^D)^2$ converges weakly-* in $L^\infty((0,T)\times(0,1))$ to $(\partial_x p)^2$.
This allows us to pass to the limit in the term $(\partial_x p^D)^2 m^D$ in \eqref{eq1D_D0_2}.
Lemma~\ref{lem:m-alg-conv} establishes the limit passage in the term $|m^D|^{2(\gamma-1)}m^D$.
Finally, the uniform bound on $\partial_x m^D$ in $L^2((0,T)\times(0,1))$ implied by Lemma~\ref{lem:energy}
establishes the convergence of the term $D^2 \partial_{xx}^2 m^D$ (in the weak formulation) to zero as $D\to 0$.
\endproof

\section{Outlook and open problems}
We conclude our paper by providing a list of interesting open problems
that will be the subject of future research.

\begin{itemize}
 \item
 Is the branch of mild solutions constructed in Theorem \ref{thm:exMawhin} global in $\lambda = c^2$?
 As noted in Remark \ref{rem:globalMildEx_1D}, in the one-dimensional setting $d=1$ this is true
 due to an a priori estimate on $\partial_x p$ in $L^\infty$ and maximum principle on $m$. 
 Is this also true for $d\geq 2$?
 \item
 What can we say about the dynamics of network formation? What is the mechanism of network growth?
 \item
 In connection with the previous point, it is important to understand how do stationary states depend
 on the initial data. In the case $d=1$, $D=0$ this was analyzed in Section \ref{sec:stability1D}.
 An analysis in the multi-dimensional setting is desirable.
 \item
 In Section \ref{sec:Dtozero1D} we carried out the limit $D\to 0$ in the one-dimensional setting,
 based on an estimate on $m$ in the $BV$-space. Can this be generalized to multiple dimensions?
 \item
 In connection with the previous point, an existence theorem for weak or strong solutions of the problem with $D=0$
 is needed.
 \item
 In Section \ref{subsec:LinStabD=0} we proved existence for and classified stationary states of the system with $D=0$ and $\gamma>1$.
 This should be completed by including the case $\gamma=1$.
\end{itemize}

\medskip

\noindent{\bf Acknowledgment.}
BP is  (partially) funded by the french "ANR blanche" project Kibord:  ANR-13-BS01-0004" and by Institut Universitaire de France. 
PM acknowledges support of the Fondation Sciences Mathematiques de Paris in form of his Excellence Chair 2011.



\begin{thebibliography}{99}


\bibitem{Brezis}
H. Brezis: {\sl Analyse fonctionelle}, Dunod, Paris, 1999. 

\bibitem{Kato} J. \`Avila and A. Ponce: \emph{Variants of Kato's inequality and removable singularities.}
Journal d'Analyse Math\'ematique 91 (2003), pp. 143--178.

\bibitem{morel} M. Bernot, V. Caselles, J.-M. Morel: \emph{Optimal Transportation Networks: Models and Theory}. LNM 1955, Springer-Verlag Berlin Heidelberg, 2009.

\bibitem{corson} F. Corson: \emph{Fluctuations and Redundancy in Optimal Transport Networks.}
Physical Review Letters, 104, 048703 (2010).

\bibitem{couder} S. Bohn, B. Andreotti, S. Douady, J. Munzinger, and Y. Couder:  \emph{Constitutive property of the local organization of leaf venation networks.} Physical Review E, 65, 061914 (2002).

\bibitem{Evans}
L. C. Evans: \emph{Partial Differential Equations.}
American Mathematical Society, Providence, Rhode Island, 1998.

\bibitem{Fonseca-Leoni}  I. Fonseca and G. Leoni: \emph{Modern methods in the calculus of variations: $L^p$ spaces.}
Springer Monographs in Mathematics. Springer, New York, 2007.

\bibitem{Guidetti} D. Guidetti: \emph{On elliptic systems in $L^1$.}
Osaka J. Math. 30 (1993), pp. 397--429.

\bibitem{Hess-Kato} P. Hess and T. Kato: \emph{On some linear and nonlinear eigenvalue problems
with an indefinite weight function.} Comm. in PDE (1980) 5(10), pp. 999--1030.

\bibitem{Hu} D. Hu: \emph{Optimization, Adaptation, and Initialization of Biological Transport Networks.}
Notes from lecture (2013).

\bibitem{Hu-Cai} D. Hu and D. Cai: \emph{Adaptation and Optimization of Biological Transport Networks.}
Phys. Rev. Lett. 111 (2013), 138701.

\bibitem{magnasco} E. Katifori, G. J. Sz\"ollosi and M. O. Magnasco: \emph{Damage and Fluctuations Induce Loops in Optimal Transport Networks.} Physical Review Letters, 104, 048704 (2010).

\bibitem{Marino} F. Marino: \emph{$L^{p,\lambda}$ regularity for divergence form elliptic equations
with discontinuous coefficents.} Le Mathematiche Vol. LVII (2002), Fasc. I, pp. 149--165.

\bibitem{Mawhin} J. Mawhin: \emph{Leray-Schauder degree: a half century of extensions and applications.}
Topological Methods in Nonlinear Analysis, Journal of the Juliusz Schauder Center,
Volume 14 (1999), pp. 195--228.

\bibitem{Meyers} N. Meyers: \emph{An $L^p$-estimate for the gradient of solutions of second order elliptic divergence equations.}
Annali della Scuola Normale Superiore di Pisa, Classe di Scienze $3^e$ s\'erie, tome 17, no. 3 (1963), pp. 189--206.

\bibitem{Rabinowitz} P. Rabinowitz: \emph{Some Global Results for Nonlinear Eigenvalue Problems.}
J. Funct. Anal. 7 (1971), pp. 487--513.

\bibitem{Rabinowitz-MPT} P. Rabinowitz: \emph{The mountain pass theorem: Theme and variations.}
In: \emph{Differential Equations}, Lecture Notes in Mathematics 957 (1982), pp. 237--271.

\bibitem{Sarason} D. Sarason: \emph{Functions of vanishing mean oscillation.}
Trans. AMS 207 (1975), pp. 391--405.

\bibitem{Simon87} J. Simon: \emph{Compact sets in the space $L^p(0,T; B)$.}
Ann. Mat. Pure Appl. IV (146), 1987, pp. 65-96.
\end{thebibliography}
\end{document}